\documentclass[a4paper,12pt]{article}

\usepackage{amscd}
\usepackage{amsmath}
\usepackage{amssymb}
\usepackage{amsthm} 
\usepackage[english]{babel}
\usepackage{calc}
\usepackage{color}
\usepackage{colortbl}
\usepackage{dsfont}
\usepackage{enumerate}
\usepackage{epsfig}
\usepackage[left=2.5cm,right=2.5cm,top=3cm,bottom=3cm]{geometry}
\usepackage{graphics}
\usepackage{graphicx}
\usepackage[latin1]{inputenc}
\usepackage{latexsym}
\usepackage{pdfpages}
\usepackage{pgfplots}
\usepackage{pstricks}
\usepackage{pxfonts}
\usepackage{tikz}

\usepackage{multirow}

\theoremstyle{definition}
\newtheorem{Definition}{Definition}[section]

\theoremstyle{remark}
\newtheorem{Remark}[Definition]{Remark}
\theoremstyle{plain}

\newtheorem{Theorem}[Definition]{Theorem}

\renewcommand{\div}{\mathrm{div}}
\newcommand{\R}{\mathbb{R}}

\renewcommand{\d}{{\,\mathrm{d}}}

\newcommand{\domain}{\Omega}

\newcommand{\materialIndex}{m}

\newcommand{\indexAffineDirection}{i}

\newcommand{\displacement}{u}

\newcommand{\displacementPeriodic}[1]{\tilde \displacement^{ {#1} } }

\newcommand{\displacementAffine}[1]{\hat \displacement^{ {#1} }}

\newcommand{\displacementAffineDerivativeDirection}[2]{A^{#1}_{#2}}

\newcommand{\characteristicFctVariable}{\chi}

\newcommand{\penaltyPerimeter}{\eta}

\newcommand{\object}{{\mathcal{O}}}

\def\cf{\emph{cf}}

\pgfplotsset{compat=newest}
\date{}
\title{Simultaneous elastic shape optimization for a domain splitting in bone tissue engineering}
\author{Patrick Dondl\footnote{Institute for Applied Math, University of Freiburg (Germany)}, Patrina S.~P.~Poh\footnote{Julius Wolff Institute, Charit{\'e} -- Universit{\"a}tsmedizin Berlin (Germany)}, Martin Rumpf and Stefan Simon\footnote{Institute for Numerical Simulation, University of Bonn (Germany)}}

\begin{document}

\maketitle

\begin{abstract}
This paper deals with the simulateneous optimization
of a subset $\object_0$ of some domain $\Omega$ and its complement $\object_1 = \Omega \setminus \overline \object_0$ both considered  
as separate elastic objects subject to a set of loading scenarios. If one asks for a configuration which minimizes the maximal elastic cost functional 
both phases compete for space since elastic shapes usually get mechanically more stable when being enlarged.
Such a problem arises in biomechanics where a bioresorbable polymer scaffold is implanted in place of lost bone tissue and in a regeneration phase
new bone tissue grows in the scaffold complement via osteogenesis. In fact, the polymer scaffold should be mechanically stable to bear loading in the early stage regeneration phase and at the same time the new bone tissue grown in the complement of this scaffold should as well bear the loading. 
Here, this optimal subdomain splitting problem with appropriate elastic cost functionals 
is introduced and existence of optimal two phase configurations is established for a regularized formulation.
Furthermore, based on a phase field approximation a finite element discretization is derived. 
Numerical experiments are presented for the design of optimal periodic scaffold microstructure.
\end{abstract}

\section{Introduction}

In this paper we investigate an elastic shape optimization with two competing objects $\object_0$ and $\object_1$ 
which result from a splitting of a given domain $\Omega\subset\R^d$, i.e. $\object_1 = \Omega \setminus \overline \object_0$, where $d=3$ 
turns out to be the interesting and application relevant case.
Both objects obey a constitutive law of linearized elasticity with in general different elasticity tensors and corresponding boundary conditions.
The stored elastic energy of an object is a measure for its elastic rigidity. Hence, for each object we take into account a cost functional
which is a monotone decreasing function of the stored elastic energies for a set of loading scenarios induced by prescribed boundary displacements.  
As the total objective functional we then consider the maximum of the resulting cost function values for the two objects. 

Typically, the optimization of a single object $\object \subset \Omega$ is complemented by a volume constraint or a volume penalty \cite{Al04}. Otherwise, usually 
$\object=\Omega$ would be maximally rigid and thus optimal. An equality constraint for the volume can either be ensured by a Cahn--Hilliard-type $H^{-1}$-gradient flow \cite{ZhWa07} or a Lagrange-multiplier ansatz \cite{AlJoTo03,BoCh03,LiKoHu05}. An inequality constraint $|\object|\leq V$ has been implemented in \cite{WaWaGu03, WaZhDi04, GuZhWa05} using a Lagrange multiplier approach.
In our setup the two objects simultaneously compete for space and we explicitly ask for an optimal distribution of space to the two objects.

Such a class of shape optimization problems is generically ill-posed and one observes the onset of microstructures along a minimizing sequence of the
the cost functional which are associated with a weak but not strong convergence of the characteristic functions of the elastic objects
along a minimizing sequence.  
To avoid this ill-posedness we replace the void by some weak material and add a perimeter penalty to the objective functional (cf. \cite{AmBu93} for a scalar problem).  
Alternatively, one might relax the problem, explicitly allowing for microstructures and considering a quasi--convexification of the integrand of the cost by taking the infimum over all possible microstructures.
For this relaxation rank-$d$ sequential laminates are known to be optimal for compliance minimization for a single load scenario \cite{AlBoFr97}.
Even in the multiple load scenario sequential laminates of possible higher rank are optimal.

To compute approximate solutions of elastic shape optimization problems there are different alternatives to an explicit discretization of the objects to be optimized.
An implicit representation of shapes via level sets  
\cite{Al04,AlBoFr97,AlGoJo04,AlJoTo02} can be used and combined with a
topological optimization \cite{AlJoMa04}. In \cite{AlDa14}, Allaire\,et\,al.  studied the optimization of multiphase materials with a 
regularization of the shape derivative in the level set context. 
An implicit description of shapes via phase fields---the approach also employed in this paper---is both analytically and numerically attractive. 
Phase-field models originated in the physical description of multiphase materials, where the integrant of a chemical bulk energy 
has two minima corresponding to the two phases and an additional diffusive interfacial energy encodes the preference for a small material interface.
This approach has been employed to elastic shape optimization by Wang and Zhou \cite{ZhWa07} and Blanck et al.
\cite{BlGa14, BlGa16}.
Guo\,et\,al.\ \cite{GuZhWa05} described the characteristic function of the object by the concatenation of a smoothed Heaviside function with a level-set function, where the smoothed Heaviside function acts like a phase-field profile. Wei and Wang \cite{WeWa09} encoded the object via a piecewise constant level set function closely related to the phase-field approach.

Our set-up is similar to the study of two-phase materials, where one phase is for example a good electrical, but poor heat conductor, and vice-versa for the other phase. The problem of finding a microstructure that maximizes the sum of both heat and electrical conductivity in this setting has been studied in \cite{Torquato:2002gh}. In particular, it has been conjectured that in this scalar case, domains bounded by periodic minimal surfaces (for example the Schwarz P surface) are optimal \cite{Torquato:2004ey}. Further analysis in \cite{Silvestre:2007hn}, however, casts doubt on this conjecture. 

Motivated by an application from biomechanics and medicine, in the present work we consider an elasticity state equation, as opposed to a scalar problem, and more general objective functions where multiple load cases can be combined.

The paper is organized as follows. In Section~\ref{sec:bio} we discuss the biomechanical application of an optimal design of a  polymer scaffolds for bone tissue engineering.  Then the associated state equations are described in Section~\ref{sec:state} and in Section~\ref{sec:cost}  we derive a suitable cost functional for the simultaneous optimization of both phases. An existence result for a regularized problem with soft instead of void material outside the actual objects is given in Section~\ref{sec:exist}. Section~\ref{sec:bone} deals with the actual application of optimal polymer scaffold microstructures and in Section~\ref{sec:fem} we derive a phase field approximation and discretize it based on a finite element approach. Finally, in Section~\ref{sec:num} numerical results are presented for different material properties of the two phases and different sets of loading scenarios and Section~\ref{sec:conclusions} provides some conclusions and mentions a number of open questions.

\section{A biomechanical optimization problem}\label{sec:bio}
As an application of simultaneous two phase optimization let us consider the optimal design of polymer scaffolds for bone tissue engineering. Globally, bone loss due to trauma, osteoporosis, or osteosarcoma comprises a major reason for disability. To this day, autograft, i.e., a graft of bone tissue from a different place in the same body, remains the gold standard for large scale bone loss.
This might be impossible for example due to donor site morbidity and limited availability.
Therefore, a number of substitutes are being explored \cite{Bhatt12, Campana14}. Among these substitutes, polymer scaffolds that function as a tissue expander (creating initial void space and allowing for tissue in-growth) show tremendous potential for bone regeneration \cite{Teo:2015fu, Goh:2015jy, Schuckert:2008hg,Schantz:2006ka}.

\begin{figure}[!htbp]
\centering{
\resizebox{0.8\textwidth}{!}{
  \begin{tabular}{ c  c  c }
               \begin{minipage}{0.25\textwidth} {\includegraphics[width=0.9\textwidth]{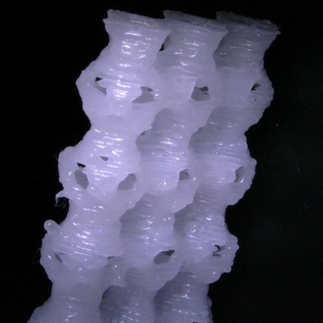}} \end{minipage}
              & \begin{minipage}{0.25\textwidth} {\includegraphics[width=0.9\textwidth]{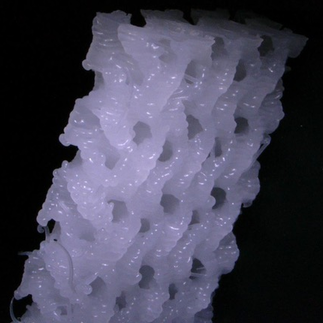}}\end{minipage}
              & \begin{minipage}{0.25\textwidth} {\includegraphics[width=0.9\textwidth]{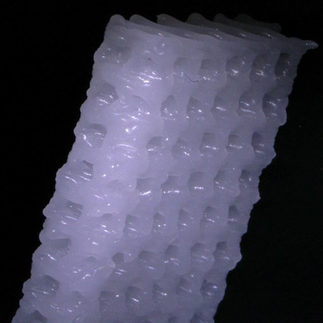}} \end{minipage}
              
  \end{tabular}
}
}
 \caption{Three different examples of 3D printed polymer scaffolds based on periodic minimal surface designs (from left to right, Schwarz P, Gyroid, and Schwarz D) with a unit cell size of approximately $1\, mm^3$. Pictures courtesy of D.~Valainis.}
 \label{fig:PrintedScaffolds}
\end{figure}

An ideal scaffold must satisfy a number of different criteria, apart from the requirement of biocompatibility.  In the initial phase of regeneration, the scaffold must provide adequate mechanical stability and the appropriate mechano-biological signal to promote osteogenesis. After new bone is formed and functional bridging through the bone defect is achieved, the scaffold should be completely resorbed allowing for a restoration of the original skeleton. With the advent of easily accessible additive manufacturing technology one possibility for such bone scaffolds are porous structures made from bioresorbable polymers, for example polycaprolactone (PCL) \cite{Poh16} (cf. the examples of printer polymer scaffolds in Figure \ref{fig:PrintedScaffolds}). Over the long regeneration time scale of more than one year, in-vivo evidence shows that such polymers degrade via bulk-erosion, that is, they lose molecular weight (and therefore mechanical stability) without a significant change in the occupied volume \cite{Lam09}, before finally being completely resorbed.

Usually, the implantation of a PCL scaffold is accompanied by a metal implant to provide further stability (see \cite{Henkel13} for illustrations of such procedures in an ovine model). It would, however, be advantageous if such metal implants were not necessary and the implanted structure, together with the regrown bone tissue, were capable of bearing the occurring mechanical loading during the regeneration time. These loading conditions depend on the stresses acting at the implant site. 
For example, if a section of the tibia is to be replaced, they would consist of unilateral compression (due to the weight of the patient) and shear (due to torsion).
Optimization procedures to design microstructures for such implants with the goal of balancing mechanical stability of the scaffold and bone tissue regeneration have been explored for example in \cite{Adachi06}. 
Furthermore, as it is suggested by the optimality conjecture for competing phases \cite{Torquato:2004ey}, designs for bone scaffolds based on periodic minimal surfaces are under consideration \cite{Kapfer:2011kz}.
Currently, the structure of the printed polymer strands still limits the achievable microstructures (cf. Figure \ref{fig:PrintedScaffolds}).

In the beginning of the regeneration process, the polymer scaffold alone has to be able to withstand the given loading conditions. 
Later in the regeneration process, the regrown bone tissue (which due to the effect of bulk erosion can only grow in the space that was \emph{not} occupied by the scaffold) has to bear this mechanical loading.
 
Thus, we are led to a shape optimization, where we simultaneously optimize the shape of the polymer and its complement which will be occupied by 
bone generated via osteogenesis. Thereby, we focus solely on this optimization problem and not on the dynamical process of the osteogenesis and the dissolving of the polymer. Furthermore, we assume that the polymer scaffold forms a spatially homogeneous microstructure and ask for the optimal shape of the polymer phase in a fundamental cell of the scaffold with affine period boundary conditions.

\section{State equation} \label{sec:state}
We consider a domain $\Omega\subset \R^d$ ($d\geq2$) with Lipschitz boundary, which is split up in two subdomains $\object_0$ and $\object_1$ 
($\overline \object_0 \cup \overline \object_1 = \overline \Omega$ and $\object_0 \cap \object_1 = \emptyset$) with corresponding characteristic functions $\chi_0$
and $\chi_1$, respectively. At first, we suppose the boundary of both domains to be Lipschitz and take into account displacement $u^m: \object_m \to \R^d$,
with a decomposition  
\begin{align}
u^m = \hat u^m + \tilde u^m\,,
\end{align}
where $\hat u^m \in H^{1,2}(\object_m)^d$ is fixed and $\tilde u^m$ lies in a closed subspace $\mathcal{V}^m$ of $H^{1,2}(\object_m)^d$.
In explicit $\hat u^m$ and the choice of $\mathcal{V}^m$ determine the encountered boundary condition, e.g. 
Dirichlet or period boundary condition (see Section~\ref{sec:bone})
on $\partial \object^m \cap \partial \Omega$  and Neumann boundary conditions on $\partial \object^m \cap \Omega$. 
Here, we assume that $\mathcal{V}^m$ is such that Korn's inequality holds for displacements $u^m\in \mathcal{V}^m$, i.e. there exists a constant 
$C_K$ such that 
\begin{align*}
\|u^m\|_{L^2(\object^m)} \leq  C_K \|\varepsilon[u^m]\|_{L^2(\object^m)}\,,
\end{align*}
where the strain tensor $\varepsilon[u]$ is given as $\frac12 (Du^T+Du)$ with $Du$ defining the Jacobian of the displacement $u$.
Now, we take into account linearized elasticity and consider the elasticity tensors $C^m = (C^m_{ijkl})_{i,j,k,l = 1,\ldots, d}$ of both subdomains ($m=0,1$).
For simplicity we assume that the two materials are isotropic
and thus determined by the Lam{\'e}-Navier parameters $\mu^\materialIndex>0$ and $\lambda^\materialIndex>0$, 
i.e. $C \varepsilon[\displacement^\materialIndex] : \varepsilon[\displacement^\materialIndex] = 2 \mu \varepsilon[\displacement^\materialIndex] : \varepsilon[\displacement^\materialIndex] + \lambda \div(\displacement^\materialIndex) \div(\displacement^\materialIndex)$.
A generalization allowing for anisotropic materials is straightforward. Then, 
the associated energy of a displacement $u\in H^{1,2}(\object_m)^d$ is given by 
\begin{align}\label{eq:energy}
E^m[\chi^m,u] = \frac12 \int_\Omega \chi^m \, C^m \varepsilon[u] : \varepsilon[u] \d x \,,
\end{align}
where $C^m A:B = \sum_{i,j,k,l} C^m_{ijkl} A_{ij} B_{kl}$ for $A,B \in \R^{d,d}$. 
Thereby, $C \varepsilon[u]$ is the stress tensor associated with the strain tensor $\varepsilon[u]$.
Using Korn's inequality it is easy to see that there exists a unique 
displacement 
$u^m$ on $\object_m$ 
which minimizes the energy $E[\chi^m,\cdot]$.  
The corresponding weak form of the Euler Lagrange equations 
\begin{align}\label{eq:state}
0 = \partial_{u}E^m[\chi^m,u^m](\phi) = \int_\Omega \chi^m \, C^m \varepsilon[u^m] : \varepsilon[\phi] \d x 
\end{align}
for all $\phi \in \mathcal{V}^m$ is the state equation for $u^m$ on $\object_m$. To express the dependence of $u^m$ on the 
shape of the subdomain $\object_m$ we also write $u^m[\chi^m]$.

\section{Cost functional} \label{sec:cost}
As already discussed in the introduction the stored elastic energy of an object is a measure for its elastic rigidity. 
We assume that a set of different boundary conditions, encoded via $\hat u^m_l$ for $l=1, \ldots, L$, 
reflects typical loading conditions of the subdomain $\object_m$. 
For the sake of simplicity we consider them to be independent of the subdomain.
Thus, we take into account a continuous function $g: \R^d \to \R$, which is supposed to be monoton decreasing in each argument 
and define 
\begin{align*}
J^m[\chi^m] := g(E^m_{1}, \ldots, E^m_{d})
\end{align*} 
as the cost associated with the set of loading conditions, where 
$E^m_{l} := E^m[\chi^m, u^m_l[\chi^m]]$ is the stored elastic energy of the equilibrium solution 
$u^m_l$ of \eqref{eq:state} corresponding to the prescribed  $\hat u^m_l$.
Aiming for an optimization of the expected value of the total energy we would have to choose
\begin{align*}
g(E_1,\ldots, E_L) = - (E_1 + \ldots + E_L)\,.
\end{align*}
In this paper we take into account 
\begin{align*}
g(E_1,\ldots, E_L) = \left(\sum_{l=1,\ldots,L} (E_l)^{-p} \right)^{\frac1p}\,.
\end{align*}
For $p\to\infty$ the resulting cost converges to the maximal inverse total energy  
$$\max_{l=1,\ldots,L} (E_l)^{-1} = (\min_{l=1,\ldots,L} E_l)^{-1}$$ 
and thus represents a worst case optimization problem where solely the loading scenario with the smallest stored elastic energy is taken into account.

Now, the objective functional associated with the domain splitting $\chi^0 =\chi$ and $\chi^1=1-\chi$ for a characteristic function $\chi$
can be defined as 
\begin{align}\label{eq:totalcost}
J[\chi] = \max \left(J^0[\chi],J^1[1-\chi]\right)\,.
\end{align}
This objective functionals reflects the competition of both subdomains aiming to increase their rigidity via domain enlargement with significant payoff in the cost $J^m[\chi^m]$. 
Due to this competition no volume constraint or penalty is needed to formulate the shape optimization problem.

\section{A hard--soft approximation and a perimeter regularization}  \label{sec:exist}
It is advantageous to approximate the characteristic function $\chi^m$ for each object $\object_m$ by $\chi^m + \delta (1-\chi^m)$ for some small constant $\delta >0$.
In explicit, instead of considering the elasticity problem on $\object_m$ we take into account an elasticity problem on the domain $\Omega$. 
To obtain a well-posed optimization problem, the void phase on the complementary set $\Omega \setminus \overline \object_m$ is replaced by a very soft phase.
Consequently we consider $u^m$, $\hat u^m$, and $\tilde u^m$ as functions in $H^{1,2}(\Omega)^d$ and $\mathcal{V}^m$ as a subspace of $H^{1,2}(\Omega)^d$. 
Then, Korn's inequality is assumed to hold on this extended function space.
Furthermore, taking into account characteristic functions $\chi\in BV(\Omega, \{0,1\})$ we add a penalty given by the perimeter of the subdomains $\penaltyPerimeter |D\chi|(\Omega)$ for some $\penaltyPerimeter>0$.
Then, we obtain the following straightforward existence theorem for a minimizer. 
\begin{Theorem}[Existence of an optimal subdomain splitting]\label{thm:exist}
For $\delta >0$, given displacements $\hat u^m_l$ with $m=0,1$ and $l=1,\ldots, L$, elastic energies
\begin{align}\label{eq:energydelta}
E^{m,\delta}[\chi,u] = \frac12 \int_\Omega (\chi + \delta (1-\chi)) \, C^m \varepsilon[u] : \varepsilon[u] \d x 
\end{align}
for $m=0,1$, and objective functional 
\begin{align}\label{eq:totalcostPer}
J^\penaltyPerimeter[\chi] = \max \left(J^0[\chi],J^1[1-\chi]\right) + \penaltyPerimeter |D\chi|(\Omega)
\end{align}
there exists a characteristic function $\chi$ which minimizes $J[\cdot]$ 
over all the admissible characteristic functions in $BV(\Omega, \{0,1\})$ 
where $\tilde u^m_l[\chi^m]$ is the unique minimizer of $E^{m,\delta}[\chi^m,\cdot + \hat u^m_l[\chi]]$ over all displacements in $H^{1,2}(\Omega)^d$
with $\chi^0=\chi$ and $\chi^1=1-\chi$.
\end{Theorem}
\begin{proof} 
For fixed $\chi \in BV(\Omega, \{0,1\})$ a direct application of the Lax-Milgram theorem combined with Korn's inequality ensures the existence of 
unique minimizers $u^0[\chi]$ and $u^1[1-\chi]$ of the energy  $E^{0,\delta}[\chi,\cdot]$ in $\mathcal{V}^0 + \hat u^0$ and $E^{1,\delta}[\chi,\cdot]$ in 
$\mathcal{V}^1 + \hat u^1$ , respectively.

Now,  we consider a minimizing sequence $(\chi_k)_{k=1,\ldots}$ of the objective functional $J^\penaltyPerimeter$ in $BV(\Omega, \{0,1\})$.
Due to the perimeter term in the objective functional the $\chi_k$ are uniformly bounded in $BV(\Omega, \{0,1\})$. Hence, there is a weakly converging 
subsequence, which we denote for simplicity again by $(\chi_k)_{k=1,\ldots}$ and which converges to some $\chi \in BV(\Omega, \{0,1\})$.
By the definition of the energies $E^{m,\delta}$  the corresponding minimizing displacements $u^m_{l,k}$ for $m=0,1$, $l=1,\ldots, L$ and $k=1,\ldots$ are
uniformly bounded in $H^{1,2}(\Omega)^d$ with uniformly bounded energies $E^{m,\delta}[\chi^m_k,u^m_{l,k}[\chi^m_k]]$. 
Thus, for $m=0,1$ and $l=1,\ldots, L$ there is subsequences, again denoted by $(u^m_{l,k})_{k=1,\ldots}$, which converge weakly in $H^{1,2}(\Omega)^d$.
Using the compact embedding of $BV(\Omega)$ in $L^1$ and Lebesgue's dominated convergence theorem, one obtains the $\Gamma$-convergence of the 
functionals $E^{m,\delta}[\chi^m_k,\cdot + \hat u^m_l[\chi_k]]$ to $E^{m,\delta}[\chi^m,\cdot + \hat u^m_l[\chi_k]]$  for $k\to \infty$ in the $H^{1,2}$-topology. 
As a direct consequence of this and the equi-coerciveness of the elastic energies,
$E^{m,\delta}[\chi^m_k,\cdot + u^m_l[\chi_k]]$ converges to $E^{m,\delta}[\chi^m,\cdot + u^m_l[\chi]]$. 
Finally, the continuity of $g$ and the $\max$ function implies the lower semi-continuity of the objective functional $J^\penaltyPerimeter$. Thus the claim holds.
\end{proof}

For a similar proof in the case of nonlinear elastic shape optimization and a phase field approach instead of an approach with characteristic functions in $BV$
we refer to \cite{PeRuWi12}.

\section{Optimal microstructured polymer scaffold} \label{sec:bone}
In the context of the biomechanical application described in Section~\ref{sec:bio} we consider a microstructured scaffold.
The spatial scale of the microstructure is thereby determined by the 3D printer technology and biological considerations, 
such as the nutrient supply via
blood vessel of a minimum thickness.
Hence, we are led to the problem of an optimal domain splitting described in Sections \ref{sec:state} and \ref{sec:cost}. 
In explicit, we consider $\Omega=[0,1]^d$ as the fundamental cell of the polymer scaffold. We take into account prescribed affine displacements $\hat u^m_l$ 
with a symmetrized strain tensor $\epsilon[u^m_l] = A_l$ with $A_l \in R^{d,d}_{sym} \cap GL(d)$ and choose the subspace 
\begin{align*}
\mathcal{V}^m =  \mathcal{V} := 
 H^1_{\#}(\domain,\R^d) = \left\{ \displacementPeriodic{} \in H^{1,2}(\domain,\R^d) \, : \, \displacementPeriodic{} \text{ periodic on } \domain \,,\,\int_\domain \displacementPeriodic{} \d x = 0 \right\}\,.
\end{align*}
Then, the elastic energy of a displacement $u^m_l =  \hat u^m_l + \tilde u^m_l$ is given by
\begin{align*}
E^m[\chi^m,u^m_l] = \frac12 \int_\Omega \chi^m \, C^m \left( A_l + \varepsilon[u]\right) : \left( A_l + \varepsilon[u)]\right) \d x \,.
\end{align*}
Let us remark that it is well-known from the theory of elastic homogenization \cite{Al02} that the 
homogenized elasticity tensor $C^\materialIndex_{*}$  of the resulting microstructure is uniquely described by
\begin{align}\label{eq:homogenzedTensor}
 C^\materialIndex_{*} B : B 
 = \min_{\displacementPeriodic{\materialIndex} \in \mathcal{V}} 
  \int_\domain \characteristicFctVariable^\materialIndex C^\materialIndex \; \left(B + \varepsilon[\displacementPeriodic{\materialIndex}]\right) \; : \; \left(B + \varepsilon[\displacementPeriodic{\materialIndex}]\right)  \d x 
\end{align}
for all  $B\in R^{d,d}_{sym} \cap GL(d)$.

\section{Phase field approximation and finite element discretization} \label{sec:fem}
A direct numerical treatment of the characteristic function $\chi$ or an explicit parametric description of the 
subdomains $\object_0$ and $\object_1$  is algorithmically quite demanding. 
Hence, we replace the characteristic function $\chi$ by a phase-field function $v:\domain \to\R$ of Modica--Mortola type. 
Then the associated phase field energy functional is given by 
\begin{equation*}
L^\varepsilon[v]:=\frac12\int_\domain \varepsilon|\nabla v|^2+\tfrac1\varepsilon\Psi(v)\,\d x\,,
\end{equation*}
where $\varepsilon$ describes the width of the diffused interface between the two subdomains (\cf\  \cite{PeRuWi12}).
Here, we set $\Psi(v):=\frac{9}{16}(v^2-1)^2$ 
with  two minima at $v=-1$ and $v=1$ and replace the perimeter $|D\chi|(\Omega)$ by the phase field energy $L^\varepsilon[v]$.
In the limit $\varepsilon\to0$, the phase field $v$ leads to a clear separation between two pure phases $-1$ and $1$ 
and $L^\varepsilon$ is $\Gamma$-converging to the perimeter functional $|D\chi|(\Omega)$ of both faces in the domain $\Omega$ \cite{Br02}.

For the numerical discretization in 3D ($d=3$) we use a cuboid mesh, i.e. the unit cube $\domain$ is uniformly divided into $(N - 1)^3$ cuboid elements with $N^3$ nodes. 
On this mesh we define the space $\mathcal{V}_h$ of piecewise trilinear, continuous functions and consider discrete phase fields $v\in  \mathcal{V}_h$ and 
discrete displacement $\displacement^\materialIndex_h \in \mathcal{V}_{h}^3$. 
In analogy to the continuous case we then restrict to space of discrete, affine periodic functions. 
Furthermore, the elastic energies are approximated by a tensor product Simpson quadrature.
Concerning the solver, the average value conditions on $\displacement^\materialIndex_h$ are imposed via a Lagrange multiplier approach. 
The corresponding linear systems for the elasticity problems are solved using the conjugate gradient method with diagonal preconditioning.
The actual shape optimization problem in the unknown phase field $v_h$ is solved using the IPOPT package \cite{WaechterBiegler2006}.
To implement the periodicity we identify the nodal values of the discrete phase field and the discrete displacements on corresponding pairs of nodes.
To deal with the translational invariance of the phase field description of the subdomains -- indeed if $v$ is optimal, then the periodically extended 
$v(\cdot + \xi)$ is also optimal for all $\xi$ -- we fix the center of mass of the phase field $v$ taking into account the additional constraints
$\int_\domain \chi \, v_h (x_i - \tfrac12) \d x = 0$ for $ i=1,2,3$. 

\section{Numerical results} \label{sec:num}
\definecolor{activeLoadsShapeOptBones}{rgb}{0.95, 0.95, 0.95}
In this section we present computational results for optimal microstructures in 3D and their dependence on material and model parameters.
Furthermore, in a conceptual study we investigate realistic material parameters for bone material and a bioresorbable polymer.

The computational results are obtained on mesh with $65^3$ or $33^3$ vertices, where we use a prolongation of the optimal phase field on $17^3$ mesh  
as the initialization. On this coarser mesh random values in the interval $[-1,1]$ are used to initialize the phase field. 
For the phase field parameter we choose $\epsilon=2h$ where 
$h$ is the grid size.

To take into account compression and shear modes in the cost functional we investigate different sets of load scenarios based on the following $6$ 
affine displacements
$\displacementAffine{}_\indexAffineDirection (x) = \displacementAffineDerivativeDirection{}{\indexAffineDirection} x$ 
as boundary data with 
$A_{11} = \beta e_1^T e_1$,  $A_{22} = \beta e_2^T e_2$, $A_{33} = \beta e_3^T e_3$, $A_{12} = \beta (e_1^T e_2 + e_2^T e_1)$, 
$A_{13} = \beta (e_1^T e_3 + e_3^T e_1)$, and $A_{23} = \beta (e_2^T e_3 + e_3^T e_2)$ for $\{e_1, e_2, e_3 \}$ being the canonical bases in $\R^3$.
We take into account the parameters $\beta = - 0.25$, $\penaltyPerimeter= 2$, and $\delta=10^{-4}$. 

Given an effective elasticity tensor $C_*$ the components $C_*^{iiii}= \beta^{-2} C_* A_{ii} : A_{ii}$ ($i=1,2,3$) represent compressive stresses caused by corresponding compressive strains,
whereas the components  $C_*^{ijij}= \beta^{-2} C_* A_{ij} : A_{ij}$ ($i,j=1,2,3,\; i\neq j$) represent shear strains induced shear stresses. 
If not indicated elsewise, we always consider $p=2$ in the definition of the weight function $g$.

\begin{figure}[!htbp]
\resizebox{1.0\textwidth}{!}{
  \begin{tabular}{ | c | c  c | c  c | c  c | }
               \hline
              & \multicolumn{2}{|c|}{ $3 \times$ compr}
              & \multicolumn{2}{|c|}{ $2 \times$ compr, $1 \times$ shear}
              & \multicolumn{2}{|c|}{ $1 \times$ compr, $2 \times$ shear} 
              \\  \hline
              & & & & & & \\
              single cell
              & \begin{minipage}{0.125\textwidth} {\includegraphics[width=0.9\textwidth]{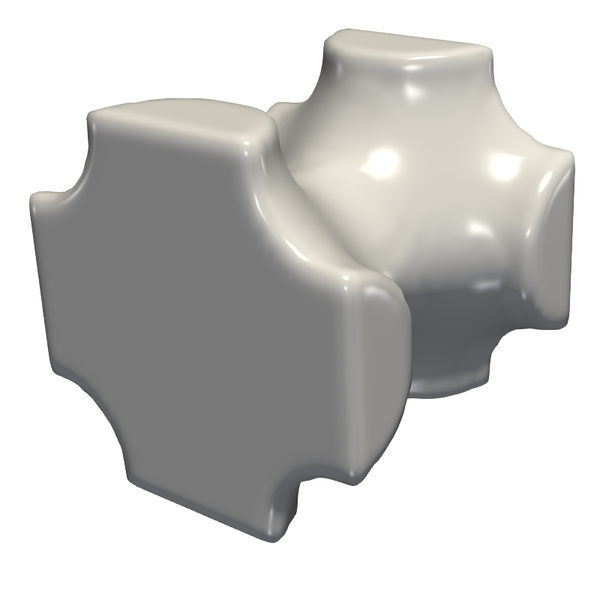}} \end{minipage}
              & \begin{minipage}{0.125\textwidth} {\includegraphics[width=0.9\textwidth]{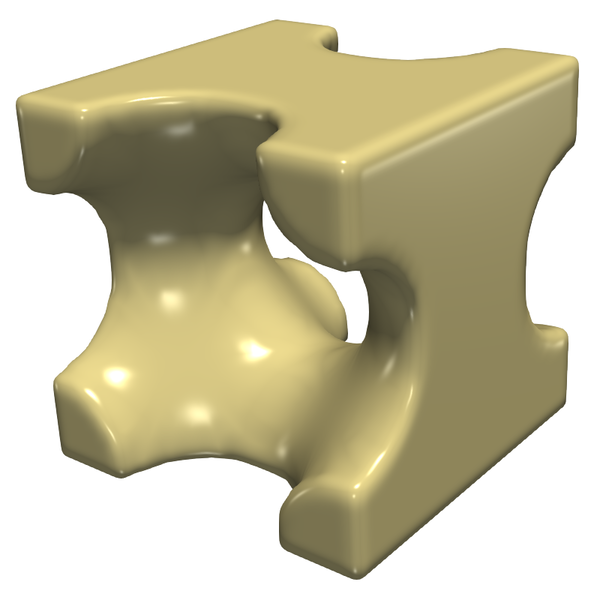}} \end{minipage}
              & \begin{minipage}{0.125\textwidth} {\includegraphics[width=0.9\textwidth]{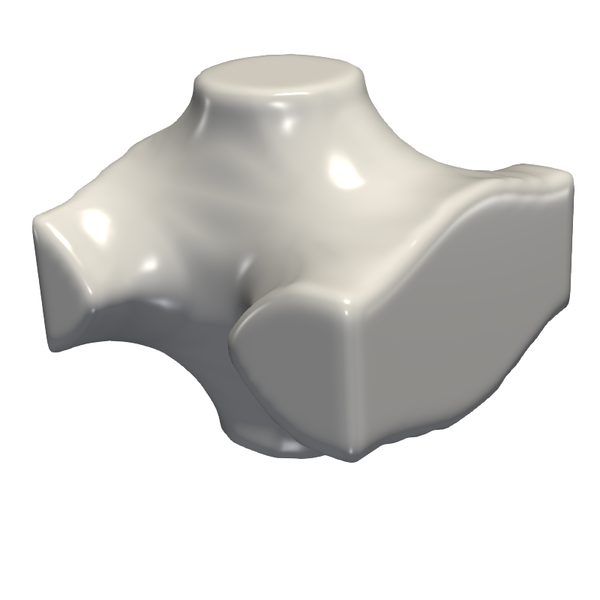}} \end{minipage}
              & \begin{minipage}{0.125\textwidth} {\includegraphics[width=0.9\textwidth]{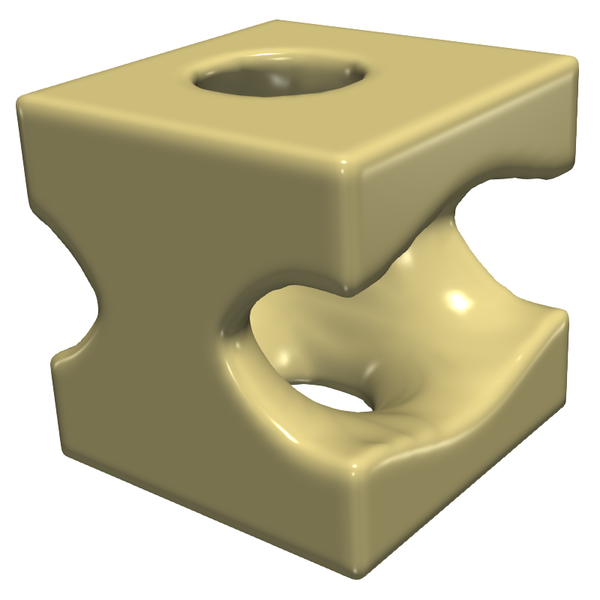}} \end{minipage}
              & \begin{minipage}{0.125\textwidth} {\includegraphics[width=0.9\textwidth]{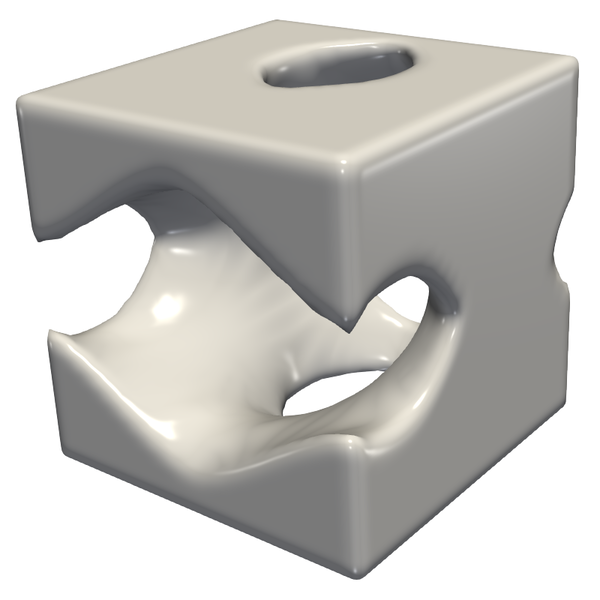}} \end{minipage}
              & \begin{minipage}{0.125\textwidth} {\includegraphics[width=0.9\textwidth]{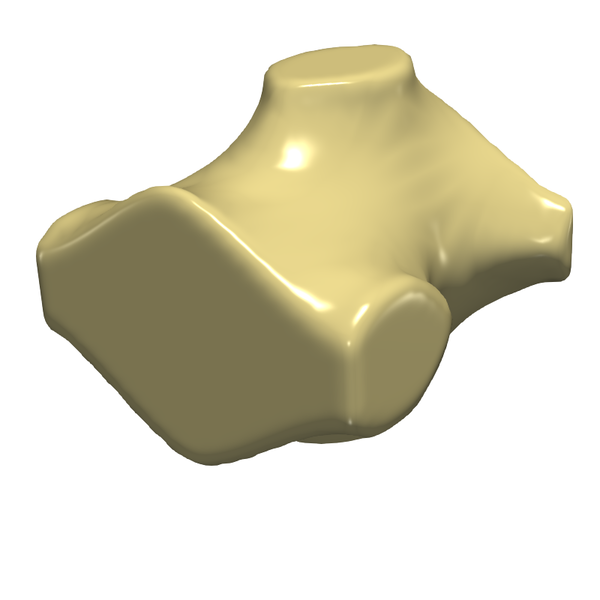}} \end{minipage}
              \\ 
             $3^3$ cells
              & \begin{minipage}{0.125\textwidth} {\includegraphics[width=0.9\textwidth]{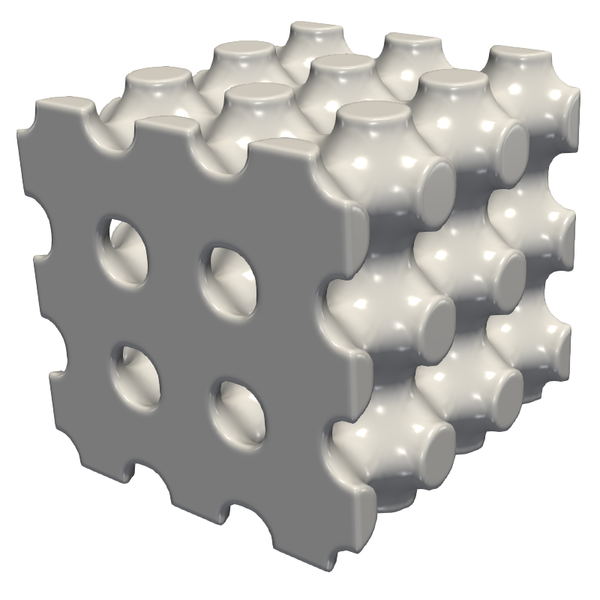}}  \end{minipage}
              & \begin{minipage}{0.125\textwidth} {\includegraphics[width=0.9\textwidth]{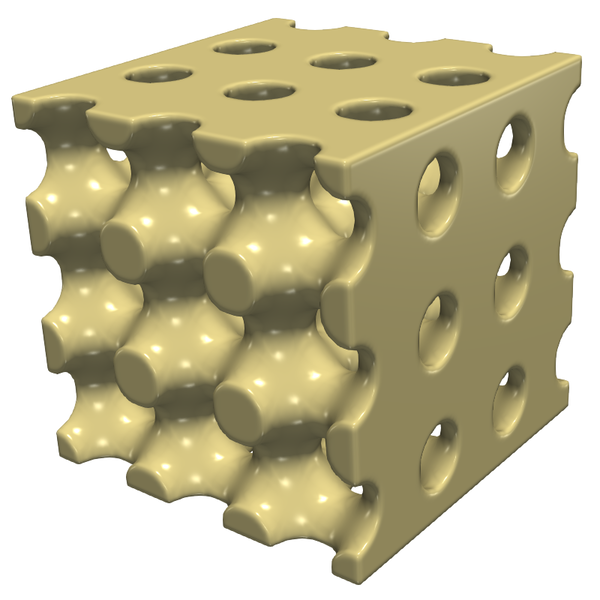}} \end{minipage}
              & \begin{minipage}{0.125\textwidth} {\includegraphics[width=0.9\textwidth]{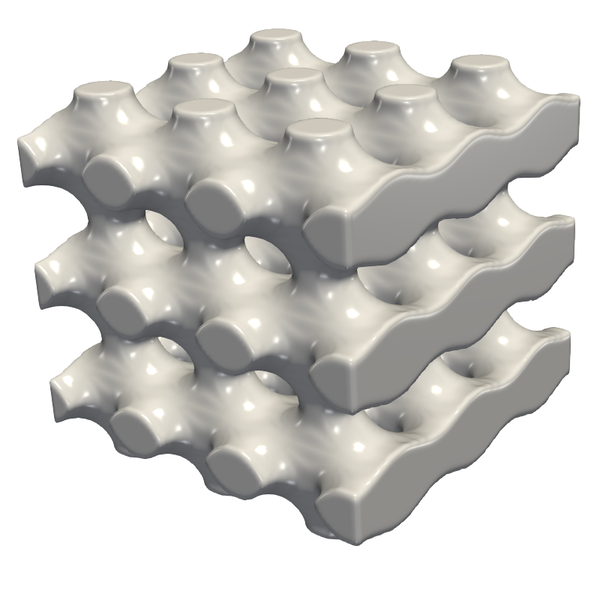}}  \end{minipage}
              & \begin{minipage}{0.125\textwidth} {\includegraphics[width=0.9\textwidth]{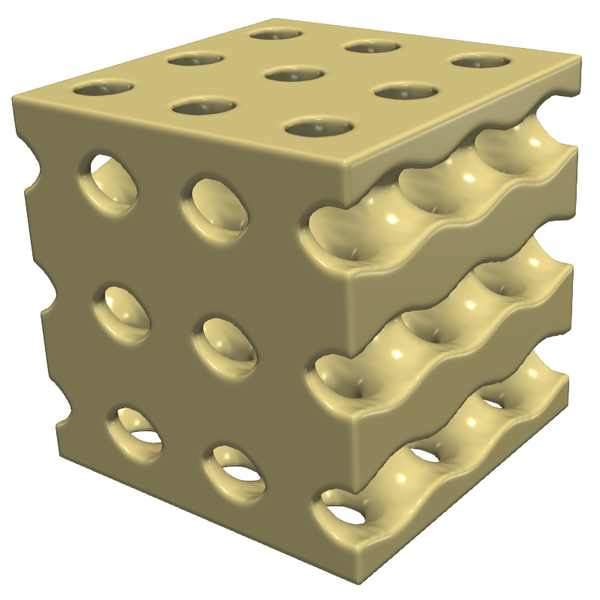}} \end{minipage}
              & \begin{minipage}{0.125\textwidth} {\includegraphics[width=0.9\textwidth]{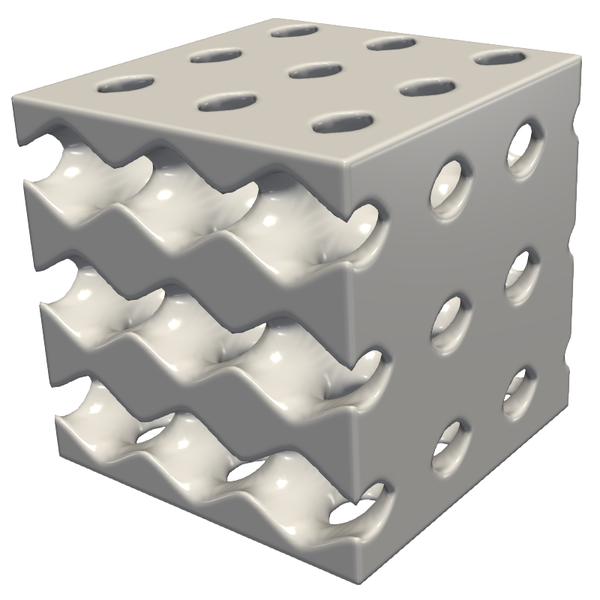}} \end{minipage}
              & \begin{minipage}{0.125\textwidth} {\includegraphics[width=0.9\textwidth]{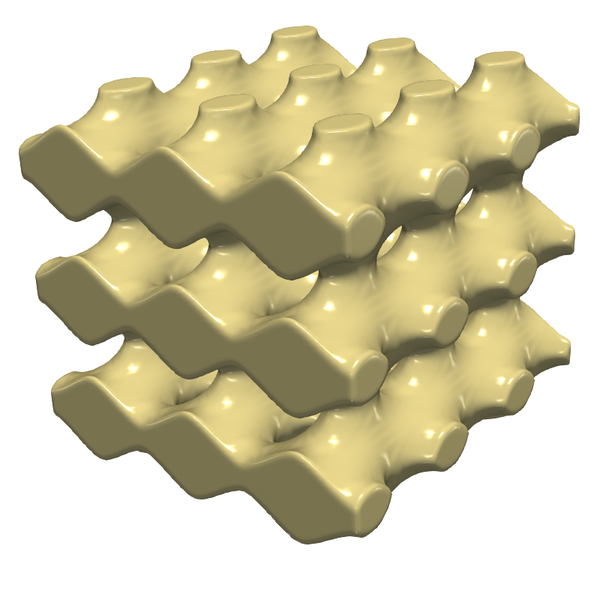}} \end{minipage}
              \\ \hline
                   & m=0 & m=1 &m=0 &m=1 &m=0 &m=1   \\ \hline
    {\small $C^{1111}_{*}$}
                   & \cellcolor{activeLoadsShapeOptBones} $ 2.825 $  & \cellcolor{activeLoadsShapeOptBones} $ 2.825 $ 
                   & \cellcolor{activeLoadsShapeOptBones} $ 2.3657 $   & \cellcolor{activeLoadsShapeOptBones}  $ 2.3657 $
                   & \cellcolor{activeLoadsShapeOptBones} $ 3.745 $  & \cellcolor{activeLoadsShapeOptBones} $ 3.745 $  \\ \hline
    {\small $C^{2222}_{*}$} 
                   & \cellcolor{activeLoadsShapeOptBones} $ 2.825 $  & \cellcolor{activeLoadsShapeOptBones} $ 2.825 $
                   & \cellcolor{activeLoadsShapeOptBones} $ 3.8584 $  & \cellcolor{activeLoadsShapeOptBones}  $ 3.8584 $
                   & $ 2.3035 $    & $ 2.3035 $ \\ \hline
     {\small $C^{3333}_{*}$} 
                   & \cellcolor{activeLoadsShapeOptBones} $ 2.825 $  & \cellcolor{activeLoadsShapeOptBones} $ 2.825 $ 
                   & $ 2.1651 $ & $ 2.1651 $
                   & $ 2.3035 $ & $ 2.3035 $\\ \hline
    {\small $C^{1212}_{*}$}
                   & $ 2.4851 $   & $ 2.4851 $ 
                   & $ 3.0126 $   & $ 3.0126 $
                   & \cellcolor{activeLoadsShapeOptBones}  $ 2.8256 $  & \cellcolor{activeLoadsShapeOptBones} $ 2.8256 $ \\ \hline
    {\small $C^{1313}_{*}$} 
                   &  $ 2.4851 $   & $ 2.4851 $ 
                   &  $ 1.1134 $   & $ 1.1134 $
                   & \cellcolor{activeLoadsShapeOptBones}  $ 2.8256 $   & \cellcolor{activeLoadsShapeOptBones} $ 2.8256 $ \\ \hline
    {\small $C^{2323}_{*}$} 
                   &  $ 2.4851 $   & $ 2.4851 $ 
                   & \cellcolor{activeLoadsShapeOptBones} $ 2.7998 $  & \cellcolor{activeLoadsShapeOptBones}  $ 2.7998 $ 
                   &  $ 1.6268 $  &  $ 1.6268 $ \\ \hline
    {\small volume}            
                   & 0.5 & 0.5 
                   & 0.5 & 0.5 
                   & 0.5 & 0.5  \\ \hline
  \end{tabular}
}
\caption{\label{fig:equalmaterial} Comparison of optimal micro-structures and relevant induced components of the effective elasticity tensors for different load scenarios indicated above.
In the top row we depict the subdomains on the fundamental cell of the microstructure and below a $3\times3\times3$ composition pronouncing the periodicity.
Those components of the tensor which are part of the corresponding objective functional are highlighted in grey.}
\end{figure}

\paragraph{Equal material parameters.}
In Fig. \ref{fig:equalmaterial} we consider equal material parameters, i.e. $(E^0,\nu^0) = (E^1,\nu^1) = (10,0.25)$, where $E^m$ is the Young's modulus 
($E^m = \frac{\mu^m ( 3 \lambda^m + 2 \mu^m}{\lambda^m + \mu^m}$) and $\nu^m$ the Poisson ratio  ($\nu^m = \frac{\lambda^m}{\lambda^m + \mu^m}$). 
Three different load scenarios are compared: three compression modes ($C_*^{1111},  \, C_*^{2222},\, C_*^{3333}$), 
two compression modes combined with a single shear mode ($C_*^{1111},  \, C_*^{2222},\, C_*^{2323}$),
and one compression mode combined with two shear modes ($C_*^{1111},  \, C_*^{1212},\, C_*^{1313}$).
For both subdomains identical loads are taken into account.
We observe an almost equal volume for both subdomains in the optimal configurations. 
In all cases the interface between the two subdomains are of the same topology as the Schwarz P surface, a periodic minimal surface representing a local minimizer of the perimeter functional. But there are significant differences in the components of the objective functional, where always those entries of the effective elasticity tensor present in the objective functional indicate a substantially stronger stiffness.

In the literature \cite{Torquato:2004ey,Silvestre:2007hn} the subdomain splitting associated with the Schwarz P surface as the interface has been 
investigated concerning its optimality in the context of a PDE constraint optimization.
On this background, we compute an optimal phase field representing a discrete minimizer of the perimeter functional as a numerical approximation of the Schwarz P surface.
For this configuration we computed the entries of the effective elasticity tensor 
and observe significantly different values $C_*^{iiii}= 2.7811$ ($i=1,2,3$) and $C_*^{ijij}= 2.481$ ($i,j=1,2,3,\; i\neq j$) compared to the optimizer in the
load scenario based on three compressions. Furthermore, the phase field  area functional $L^\varepsilon$ differs by approximately $3\%$.

\begin{figure}[!htbp]
\centering{
\resizebox{0.6\textwidth}{!}{
  \begin{tabular}{ c  c  c  c  c  c }
              \multicolumn{2}{c}{
                \begin{minipage}{0.25\textwidth}  {\includegraphics[width=0.9\textwidth]{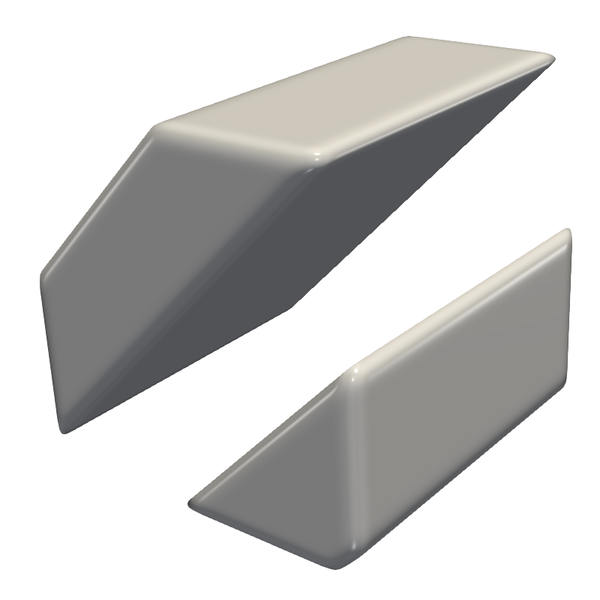}} \end{minipage}
              }
              & \multicolumn{2}{c}{
                \begin{minipage}{0.25\textwidth}  {\includegraphics[width=0.9\textwidth]{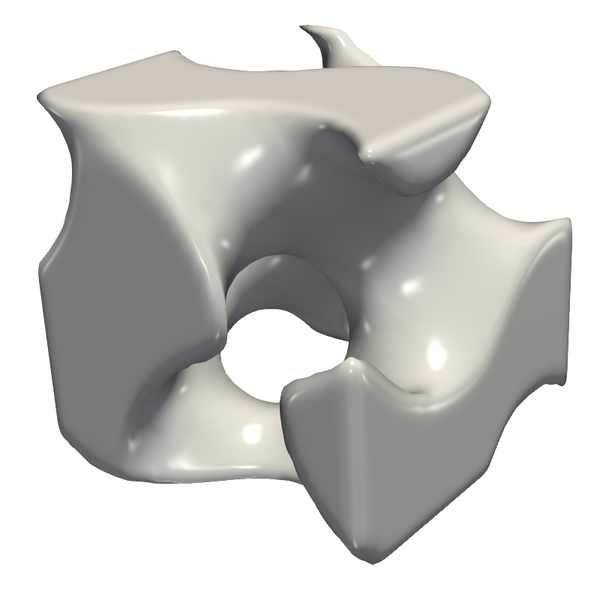}} \end{minipage}
              }
              & \multicolumn{2}{c}{
                \begin{minipage}{0.25\textwidth} {\includegraphics[width=0.9\textwidth]{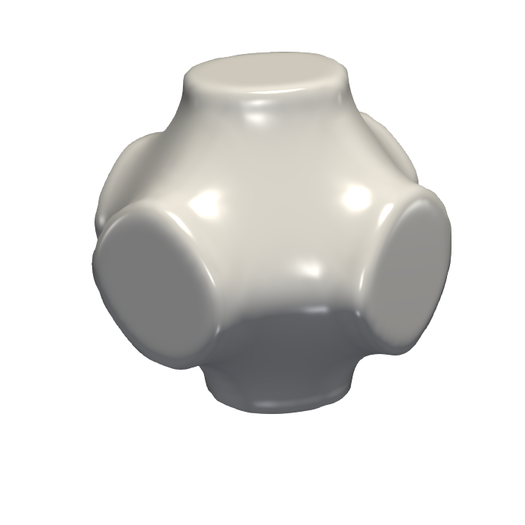}}\end{minipage}
              }  \\ 
              \multicolumn{2}{c}{
                \begin{minipage}{0.25\textwidth}  {\includegraphics[width=0.9\textwidth]{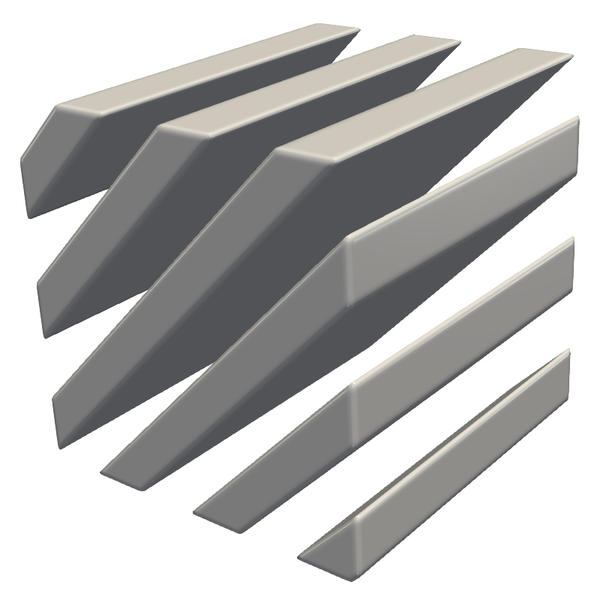}}  \end{minipage}
              }
              & \multicolumn{2}{c}{
                \begin{minipage}{0.25\textwidth} {\includegraphics[width=0.9\textwidth]{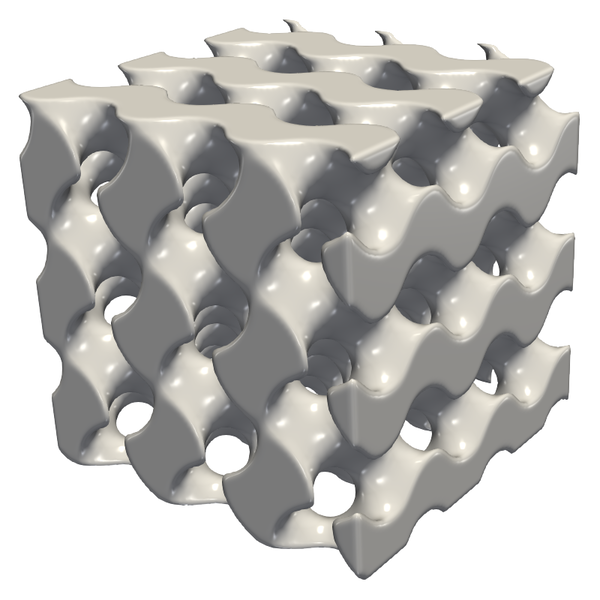}}  \end{minipage}
              }
              & \multicolumn{2}{c}{
                \begin{minipage}{0.25\textwidth}  {\includegraphics[width=0.9\textwidth]{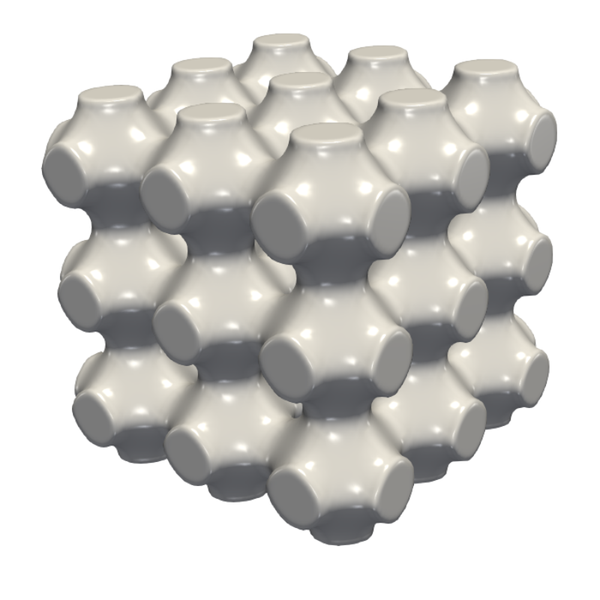}} \end{minipage}
              }  \\
  \end{tabular}
}
}
 \caption{Optimal microstructures for different values of  the perimeter parameter $\penaltyPerimeter$  (From left to right: $
  \penaltyPerimeter = 2, 4, 10$).}
 \label{fig:ShapeOptBonesPerimeter}
\end{figure}

Next, for the scenario with three shear loads ($C_*^{1212},\, C_*^{1313},\, C_*^{2323}$) 
we successively increase the parameter $\penaltyPerimeter$ in front of the perimeter functional ($\penaltyPerimeter =2,\,4,\,10$).
For small $\penaltyPerimeter$ we obtain a laminate type optimal configuration, whereas for larger $\penaltyPerimeter$ the interface is again similar to the 
Schwarz P surface as shown in Figure \ref{fig:ShapeOptBonesPerimeter}.
On the intermediate range of the parameter  $\penaltyPerimeter$ we obtain a optimal microstructure with an 
interface similar to a gyroid minimal surface \cite{Kapfer:2011kz}.

\begin{figure}[!htbp]
\resizebox{1.0\textwidth}{!}{
  \begin{tabular}{ c  c  c  c  c  c  c  c  c  c  c  }
              & \multicolumn{2}{c}{
                \begin{minipage}{0.25\textwidth}
                  {\includegraphics[width=0.9\textwidth]{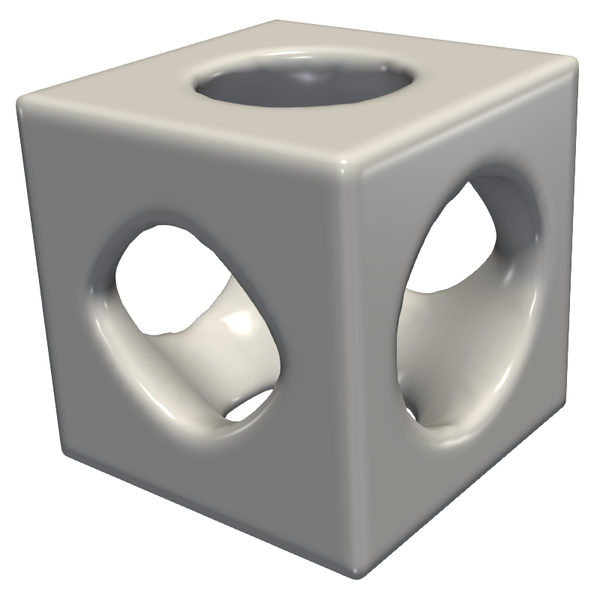}}
                \end{minipage}
              }
              & \multicolumn{2}{c}{
                \begin{minipage}{0.25\textwidth}
                  {\includegraphics[width=0.9\textwidth]{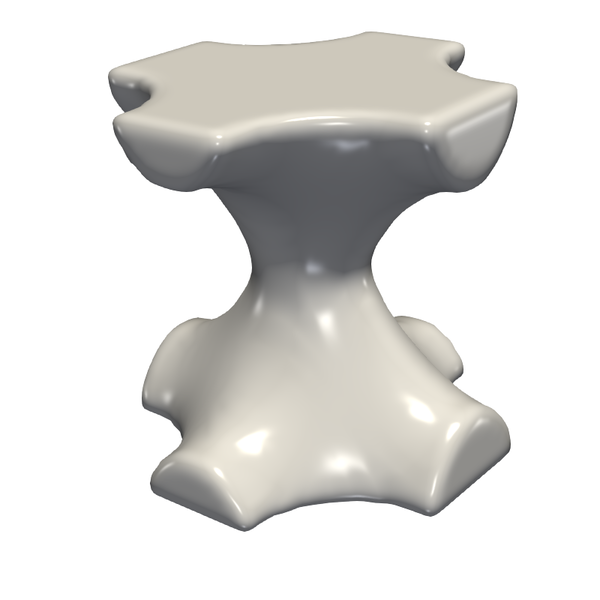}}
                \end{minipage}
              }
              & \multicolumn{2}{c}{
                \begin{minipage}{0.25\textwidth}
                  {\includegraphics[width=0.9\textwidth]{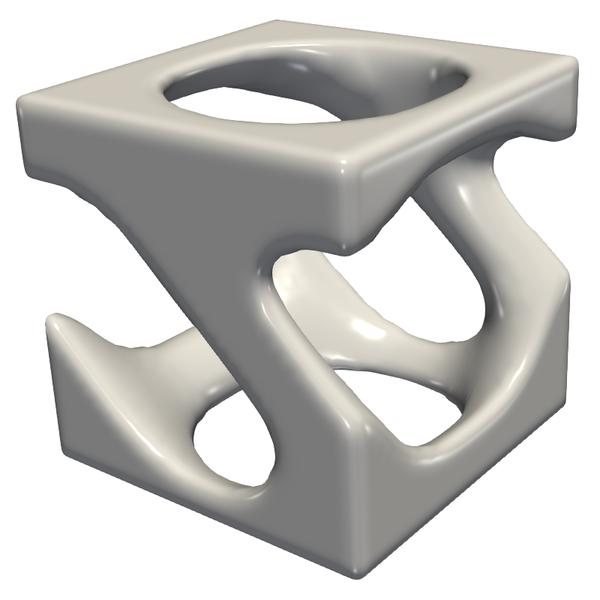}}
                \end{minipage}
              }
              & \multicolumn{2}{c}{
                \begin{minipage}{0.25\textwidth}
                  {\includegraphics[width=0.9\textwidth]{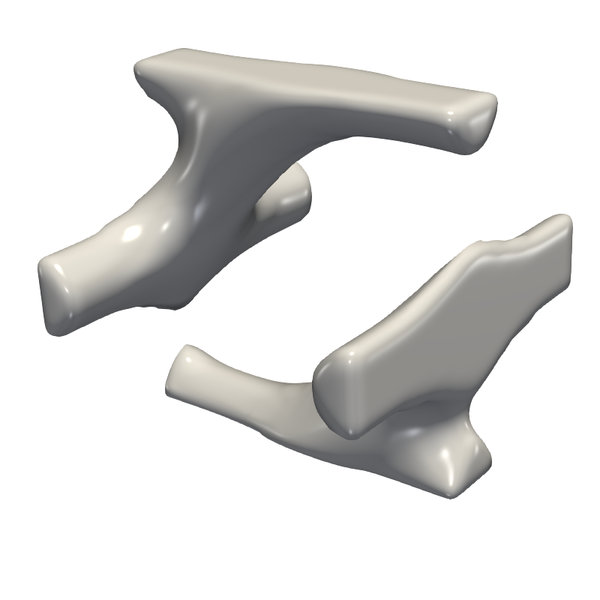}}
                \end{minipage}
              }
              & \multicolumn{2}{c}{
                \begin{minipage}{0.25\textwidth}
                  {\includegraphics[width=0.9\textwidth]{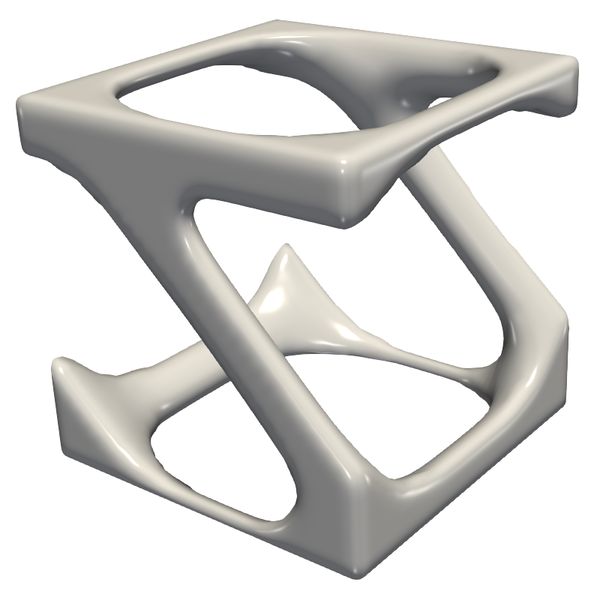}}
                \end{minipage}
              }\\ 
              & \multicolumn{2}{c}{
                \begin{minipage}{0.25\textwidth}
                  {\includegraphics[width=0.9\textwidth]{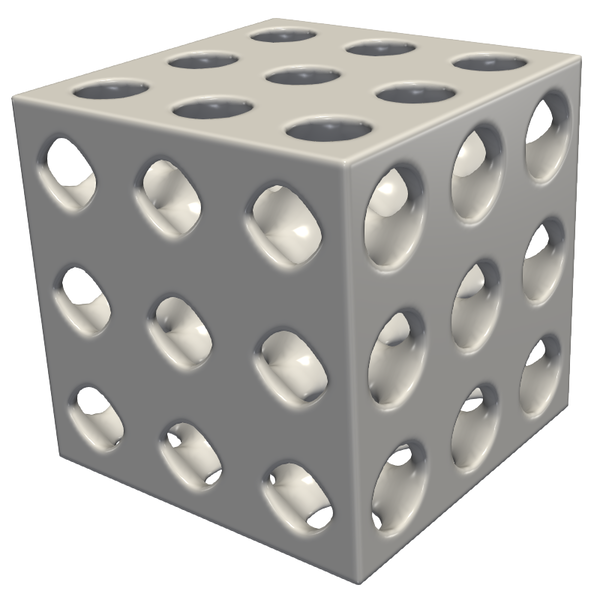}}
                \end{minipage}
              }
              & \multicolumn{2}{c}{
                \begin{minipage}{0.25\textwidth}
                  {\includegraphics[width=0.9\textwidth]{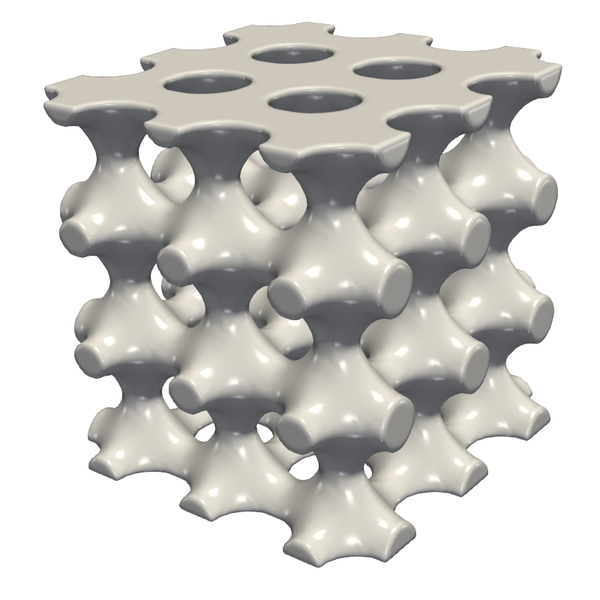}}
                \end{minipage}
              }
              & \multicolumn{2}{c}{
                \begin{minipage}{0.25\textwidth}
                  {\includegraphics[width=0.9\textwidth]{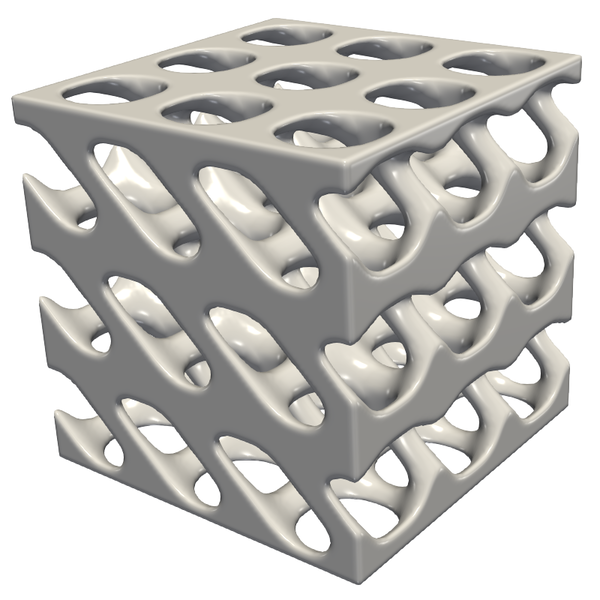}}
                \end{minipage}
              } 
              & \multicolumn{2}{c}{
                \begin{minipage}{0.25\textwidth}
                  {\includegraphics[width=0.9\textwidth]{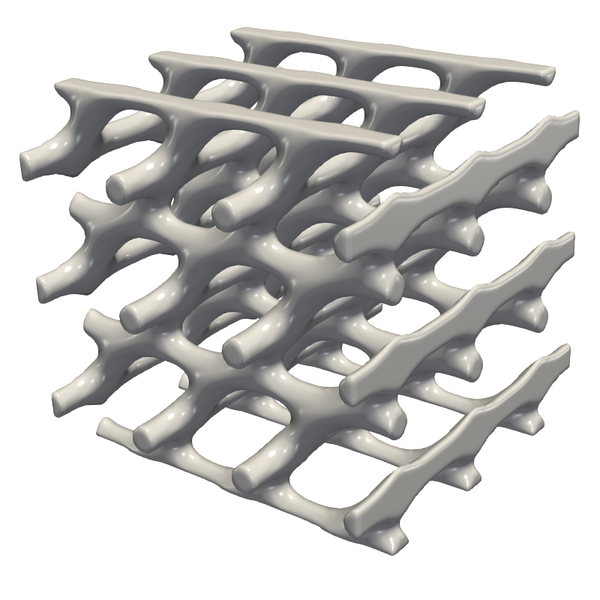}}
                \end{minipage}
              }
              & \multicolumn{2}{c}{
                \begin{minipage}{0.25\textwidth}
                  {\includegraphics[width=0.9\textwidth]{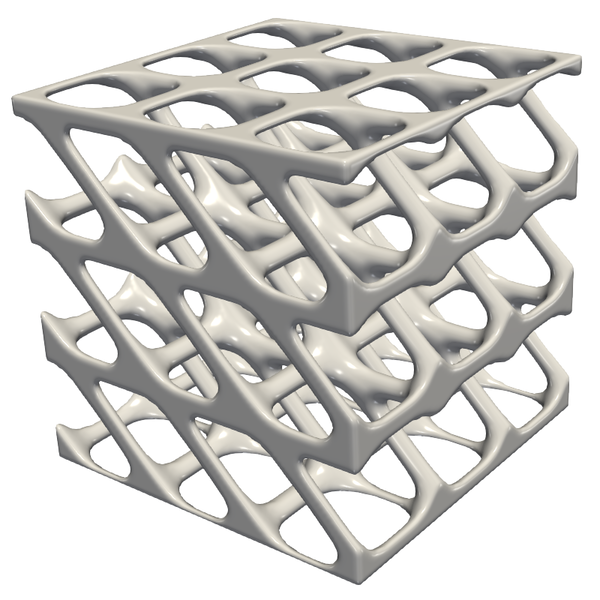}}
                \end{minipage}
              }\\ \hline
            & m=0 & m=1 &m=0 &m=1 &m=0 &m=1   &m=0 &m=1 &m=0 &m=1 \\ \hline
    $C^{1111}_{*}$ 
             & \cellcolor{activeLoadsShapeOptBones}$ 3.6587 $    & \cellcolor{activeLoadsShapeOptBones}$ 3.3904 $ 
             & \cellcolor{activeLoadsShapeOptBones}$ 4.7996 $    & \cellcolor{activeLoadsShapeOptBones}$ 4.2686 $ 
             & \cellcolor{activeLoadsShapeOptBones}$ 6.413 $    & \cellcolor{activeLoadsShapeOptBones}$ 5.9585 $ 
             & \cellcolor{activeLoadsShapeOptBones}$ 6.837 $    & \cellcolor{activeLoadsShapeOptBones}$ 7.1062 $ 
             & \cellcolor{activeLoadsShapeOptBones}$ 8.8662 $    & \cellcolor{activeLoadsShapeOptBones}$ 8.6567 $ \\ \hline
    $C^{2222}_{*}$ 
             & \cellcolor{activeLoadsShapeOptBones}$ 4.9667 $    & \cellcolor{activeLoadsShapeOptBones}$ 4.3946 $ 
             & \cellcolor{activeLoadsShapeOptBones}$ 7.0545 $    & \cellcolor{activeLoadsShapeOptBones}$ 5.4043 $ 
             & \cellcolor{activeLoadsShapeOptBones}$ 7.6653 $    & \cellcolor{activeLoadsShapeOptBones}$ 6.7548 $ 
             & \cellcolor{activeLoadsShapeOptBones}$ 7.2202 $    & \cellcolor{activeLoadsShapeOptBones}$ 7.4892 $ 
             & \cellcolor{activeLoadsShapeOptBones}$ 9.3866 $    & \cellcolor{activeLoadsShapeOptBones}$ 9.0007 $ \\ \hline
     $C^{3333}_{*}$   
             & $ 3.5702 $    & $ 3.4632 $  
             & $ 5.0673 $    & $ 4.6461 $ 
             & $ 2.3783 $    & $ 5.9513 $ 
             & $ 3.9501 $    & $ 7.6961 $
             & $ 2.562 $    & $ 8.4978 $\\ \hline
     $C^{1212}_{*}$    
             & $ 3.2681 $    & $ 4.0275 $ 
             & $ 3.3935 $    & $ 5.4312 $ 
             & $ 2.7089 $    & $ 8.1685 $
             & $ 0.67487 $    & $ 9.0034 $ 
             & $ 1.2522 $    & $ 10.1784 $ \\ \hline
    $C^{1313}_{*}$ 
             & $ 2.2299 $    & $ 3.1636 $ 
             & $ 2.3708 $    & $ 4.6865 $
             & $ 1.006 $    & $ 6.5492 $
             & $ 1.1444 $    & $ 9.4359 $ 
             & $ 0.46034 $    & $ 10.0924 $ \\ \hline
    $C^{2323}_{*}$ 
             & \cellcolor{activeLoadsShapeOptBones}$ 3.5124 $    & \cellcolor{activeLoadsShapeOptBones}$ 4.1468 $ 
             & \cellcolor{activeLoadsShapeOptBones}$ 4.283 $    & \cellcolor{activeLoadsShapeOptBones}$ 5.8455 $
             & \cellcolor{activeLoadsShapeOptBones}$ 6.0737 $    & \cellcolor{activeLoadsShapeOptBones}$ 8.1456 $
             & \cellcolor{activeLoadsShapeOptBones}$ 9.7432 $    & \cellcolor{activeLoadsShapeOptBones}$ 10.0662 $ 
             & \cellcolor{activeLoadsShapeOptBones}$ 8.8804 $    & \cellcolor{activeLoadsShapeOptBones}$ 10.1808 $ \\ \hline
    vol      & 0.41037 & 0.58963 
             & 0.32744 & 0.67256
             & 0.22565 & 0.77435
             & 0.15367 & 0.84633
             & 0.10194 & 0.89806 \\ \hline
  \end{tabular}
}
\medskip

\resizebox{1.0\textwidth}{!}{
  \begin{tabular}{ c  c  c  c  c  c  c  c  c  c  c  }
              & \multicolumn{2}{c}{
                \begin{minipage}{0.25\textwidth}
                  {\includegraphics[width=0.9\textwidth]{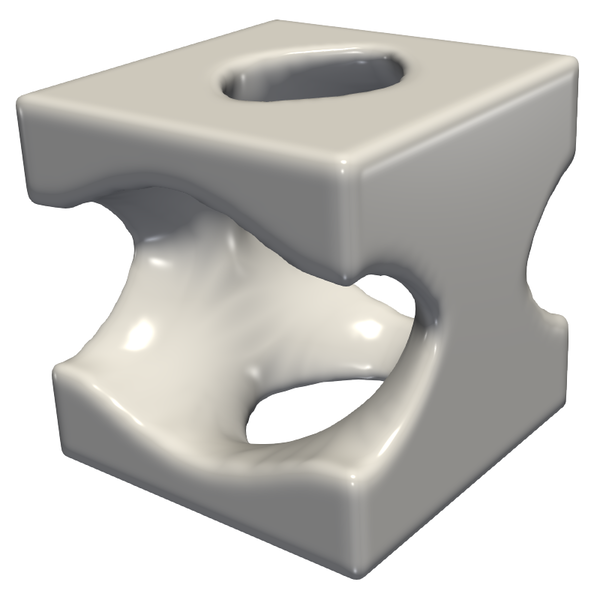}}
                \end{minipage}
              }
              & \multicolumn{2}{c}{
                \begin{minipage}{0.25\textwidth}
                  {\includegraphics[width=0.9\textwidth]{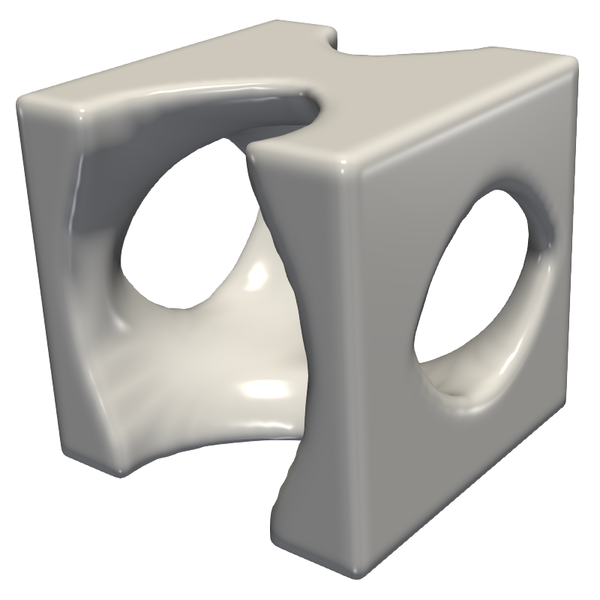}}
                \end{minipage}
              }
              & \multicolumn{2}{c}{
                \begin{minipage}{0.25\textwidth}
                  {\includegraphics[width=0.9\textwidth]{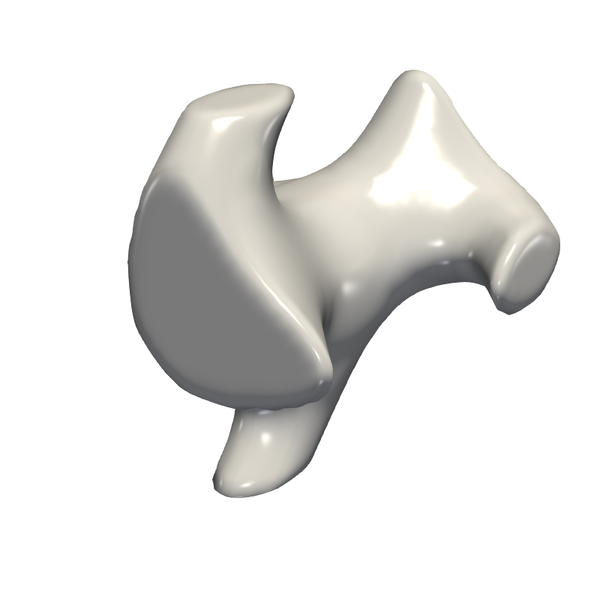}}
                \end{minipage}
              }
              & \multicolumn{2}{c}{
                \begin{minipage}{0.25\textwidth}
                  {\includegraphics[width=0.9\textwidth]{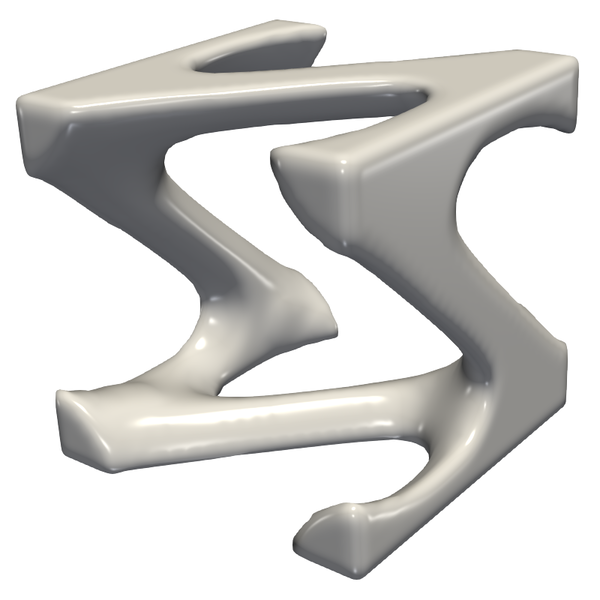}}
                \end{minipage}
              }
              & \multicolumn{2}{c}{
                \begin{minipage}{0.25\textwidth}
                  {\includegraphics[width=0.9\textwidth]{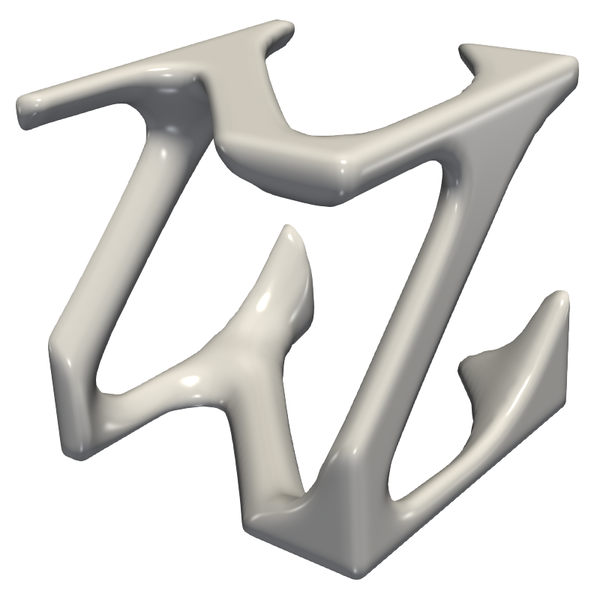}}
                \end{minipage}
              }\\ 
              & \multicolumn{2}{c}{
                \begin{minipage}{0.25\textwidth}
                  {\includegraphics[width=0.9\textwidth]{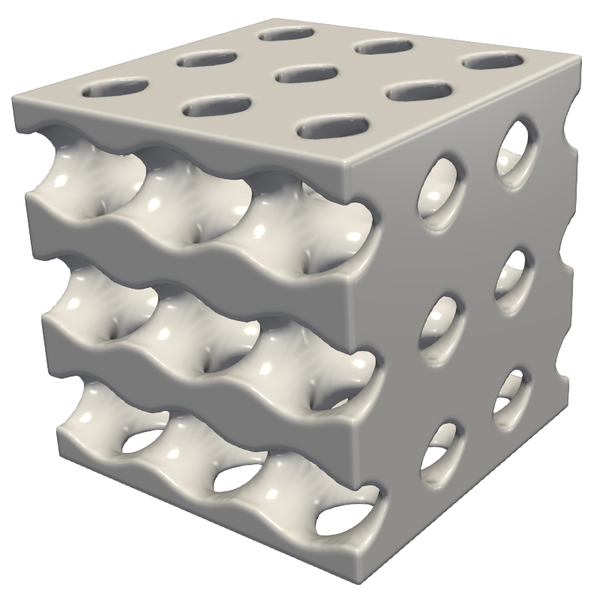}}
                \end{minipage}
              }
              & \multicolumn{2}{c}{
                \begin{minipage}{0.25\textwidth}
                  {\includegraphics[width=0.9\textwidth]{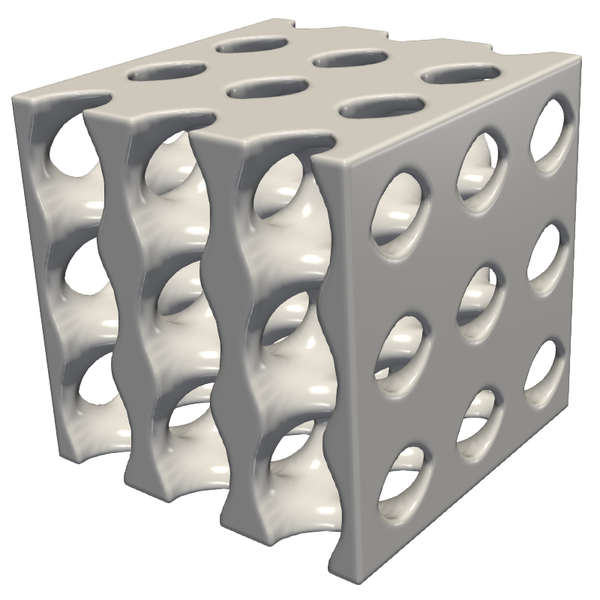}}
                \end{minipage}
              }
              & \multicolumn{2}{c}{
                \begin{minipage}{0.25\textwidth}
                  {\includegraphics[width=0.9\textwidth]{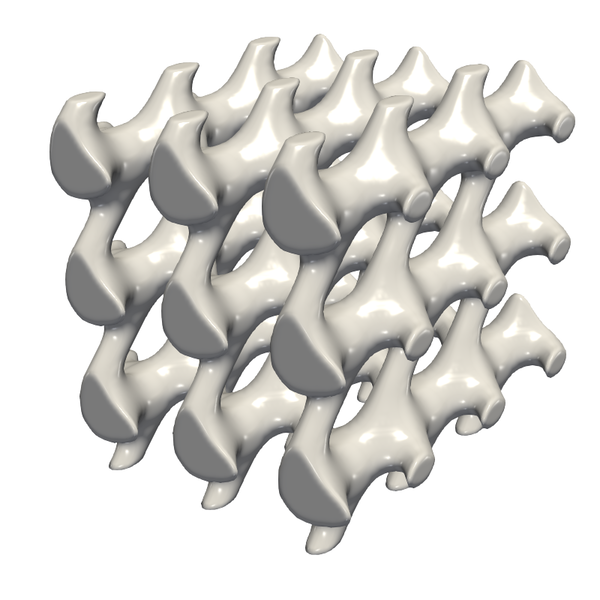}}
                \end{minipage}
              } 
              & \multicolumn{2}{c}{
                \begin{minipage}{0.25\textwidth}
                  {\includegraphics[width=0.9\textwidth]{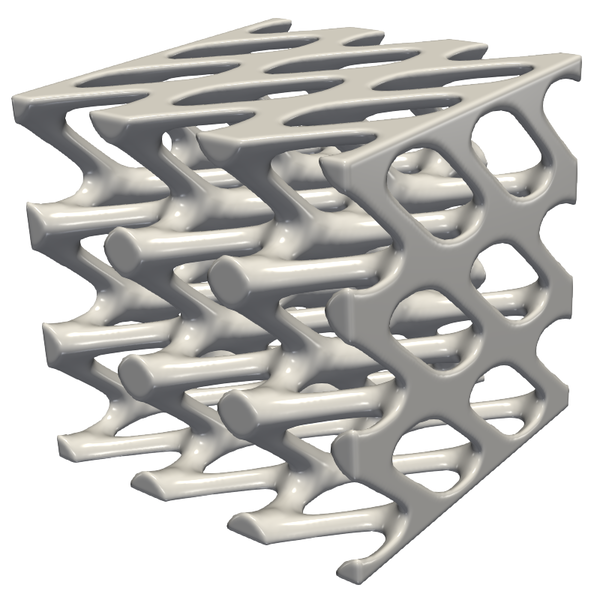}}
                \end{minipage}
              }
              & \multicolumn{2}{c}{
                \begin{minipage}{0.25\textwidth}
                  {\includegraphics[width=0.9\textwidth]{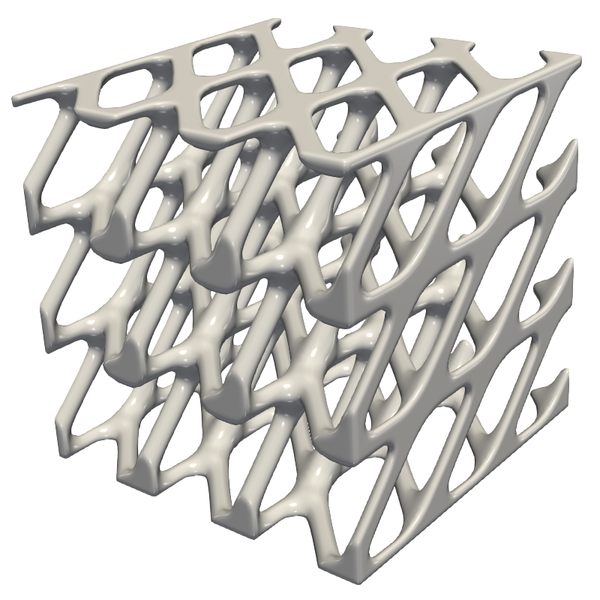}}
                \end{minipage}
              }\\ \hline
           & m=0 & m=1 &m=0 &m=1 &m=0 &m=1 &m=0 &m=1 &m=0 &m=1  \\ \hline
     $C^{1111}_{*}$ 
              & \cellcolor{activeLoadsShapeOptBones} $ 6.3411 $    & \cellcolor{activeLoadsShapeOptBones} $ 4.9912 $ 
              & \cellcolor{activeLoadsShapeOptBones} $ 9.6988 $    & \cellcolor{activeLoadsShapeOptBones} $ 5.8266 $
              & \cellcolor{activeLoadsShapeOptBones} $ 10.39 $     & \cellcolor{activeLoadsShapeOptBones} $ 6.7704 $
              & \cellcolor{activeLoadsShapeOptBones} $ 10.264 $    & \cellcolor{activeLoadsShapeOptBones} $ 8.0092 $  
              & \cellcolor{activeLoadsShapeOptBones} $ 10.773 $    & \cellcolor{activeLoadsShapeOptBones} $ 9.1178 $ \\ \hline
     $C^{2222}_{*}$ 
              & $ 3.103 $     & $ 2.8423 $  
              & $ 4.3565 $    & $ 3.6762 $ 
              & $ 3.2102 $    & $ 5.2916 $
              & $ 2.4592 $    & $ 7.0516 $ 
              & $ 2.6306 $    & $ 8.3817 $ \\ \hline
     $C^{3333}_{*}$ 
              & $ 3.1044 $    & $ 2.8416 $
              & $ 4.4754 $    & $ 3.6902 $
              & $ 3.2021 $    & $ 5.2907 $ 
              & $ 2.4592 $    & $ 7.0516 $ 
              & $ 2.6306 $    & $ 8.3817 $ \\ \hline
     $C^{1212}_{*}$ 
              & \cellcolor{activeLoadsShapeOptBones} $ 3.9041 $    & \cellcolor{activeLoadsShapeOptBones} $ 4.1364 $
              & \cellcolor{activeLoadsShapeOptBones} $ 4.8936 $    & \cellcolor{activeLoadsShapeOptBones} $ 5.4934 $  
              & \cellcolor{activeLoadsShapeOptBones} $ 6.4294 $    & \cellcolor{activeLoadsShapeOptBones} $ 7.7477 $ 
              & \cellcolor{activeLoadsShapeOptBones} $ 8.0385 $    & \cellcolor{activeLoadsShapeOptBones} $ 10.026 $ 
              & \cellcolor{activeLoadsShapeOptBones} $ 9.2991 $    & \cellcolor{activeLoadsShapeOptBones} $ 11.802 $ \\ \hline
    $C^{1313}_{*}$ 
              & \cellcolor{activeLoadsShapeOptBones} $ 3.9046 $    & \cellcolor{activeLoadsShapeOptBones} $ 4.1362 $ 
              & \cellcolor{activeLoadsShapeOptBones} $ 4.8983 $    & \cellcolor{activeLoadsShapeOptBones} $ 5.5215 $ 
              & \cellcolor{activeLoadsShapeOptBones} $ 6.4298 $    & \cellcolor{activeLoadsShapeOptBones} $ 7.7457 $ 
              & \cellcolor{activeLoadsShapeOptBones} $ 8.0385 $    & \cellcolor{activeLoadsShapeOptBones} $ 10.026 $ 
              & \cellcolor{activeLoadsShapeOptBones} $ 9.2991 $    & \cellcolor{activeLoadsShapeOptBones} $ 11.802 $\\ \hline
     $C^{2323}_{*}$
              & $ 1.2109 $    & $ 1.7901 $ 
              & $ 1.2961 $    & $ 2.7966 $ 
              & $ 1.0188 $    & $ 5.1458 $ 
              & $ 0.5915 $    & $ 8.1734 $ 
              & $ 0.36466 $   & $ 10.502 $\\ \hline
    vol       & 0.41917 & 0.58083 
              & 0.3454 & 0.6546  
              & 0.25082 & 0.74918 
              & 0.16297 & 0.83703 
              & 0.10518 & 0.89482 \\ \hline
  \end{tabular}
  }
 \caption{ \label{fig:young}
 Comparison of optimal micro-structures for varying Young's modulus  (from left to right $E^0=20,40,80,160,320$). We take into account load configurations with 
 two compression loads and one shear load (top) and  one compression load and two shear loads (bottom).
 We depict the subdomain $\object_0$ on the fundamental cell of the microstructure and a $3\times3\times3$ composition. 
 Those components of the tensor which are part of the corresponding objective functional are again highlighted in grey.}
\end{figure}

Furthermore, we investigate the impact of the choice of the weight function $g$ on the optimal splitting and the associated stiffness moduli.
In Figure \ref{fig:p} we show in the load scenario with two compression loads and one shear load the relevant entries of the effective elasticity tensor.
For increasing $p$ we observe a successive balancing of the different components of the objective functional, 
in particular the largest component $C^{2222}_{*}$ of the effective elasticity tensor is slightly deccreasing while the smallest compenent $C^{1111}_{*}$ is slightly increasing.
\begin{figure}[!htbp]
\resizebox{1.0\textwidth}{!}{
  \begin{tabular}{ c  c  c  c  c  c  c  c  c }
     \hline
        p     & \multicolumn{2}{c}{ $2$ }
              & \multicolumn{2}{c}{ $4$ }
              & \multicolumn{2}{c}{ $8$ }
              & \multicolumn{2}{c}{ $16$ } \\ \hline
        & m=0 & m=1 &m=0 &m=1 &m=0 &m=1   &m=0 &m=1 \\ \hline
   $C^{1111}_{*}$ & $ 2.3657 $    & $ 2.3657 $ 
                  & $ 2.4438 $    & $ 2.4384 $ 
                  & $ 2.4847 $    & $ 2.4808 $   
                  & $ 2.5053 $    & $ 2.5056 $ \\ \hline
   $C^{2222}_{*}$ & $ 3.8584 $    & $ 3.8584 $
                  & $ 3.8408 $    & $ 3.8429 $
                  & $ 3.8286 $    & $ 3.8291 $
                  & $ 3.8286 $    & $ 3.828 $ \\ \hline
   $C^{2323}_{*}$ & $ 2.7998 $    & $ 2.7998 $
                  & $ 2.6764 $    & $ 2.6857 $
                  & $ 2.6139 $    & $ 2.6206 $ 
                  & $ 2.5768 $    & $ 2.5766 $ \\ \hline
  \end{tabular}
}
 \caption{\label{fig:p} For different values of $p$ stiffness moduli of the optimal subdomain splitting are depicted.}
\end{figure}
\paragraph{Varying Young modulus.}
Next we study the influence of the Youngs modulus and consider $E^0 = 20,40,80,160,320$, whereas $(E^1,\nu^1) = (10,0.25)$ and $\nu^0 = 0.25$.
Furthermore,  for the perimeter parameter we choose  $\penaltyPerimeter = 1$.
We observe that the subdomain with increasing Young modulus  is getting successively thinner in the optimal domain splitting and the relative decrease in stiffness of the 
other material has to compensated by a higher volume fraction.
Figure \ref{fig:young} shows results obtained for different load scenarios.

\begin{figure}[!htbp]
\resizebox{1.0\textwidth}{!}{
  \begin{tabular}{ c  c }
       bone & polymer
       \\
          \begin{minipage}{0.25\textwidth} {\includegraphics[width=0.9\textwidth]{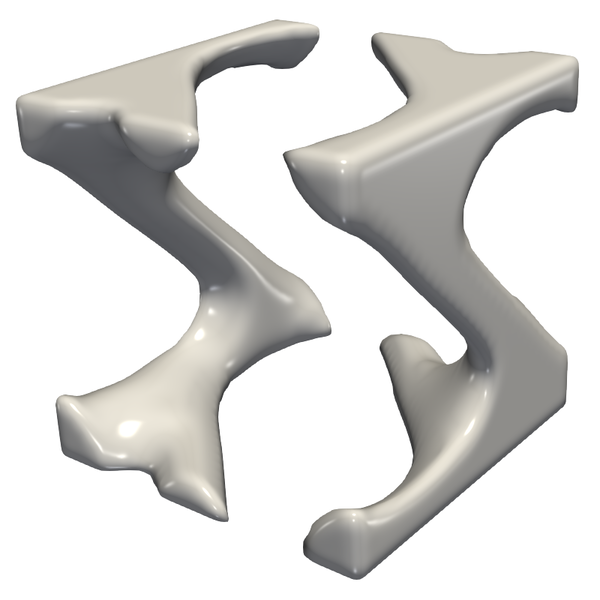}}  \end{minipage}
       &  \begin{minipage}{0.25\textwidth} {\includegraphics[width=0.9\textwidth]{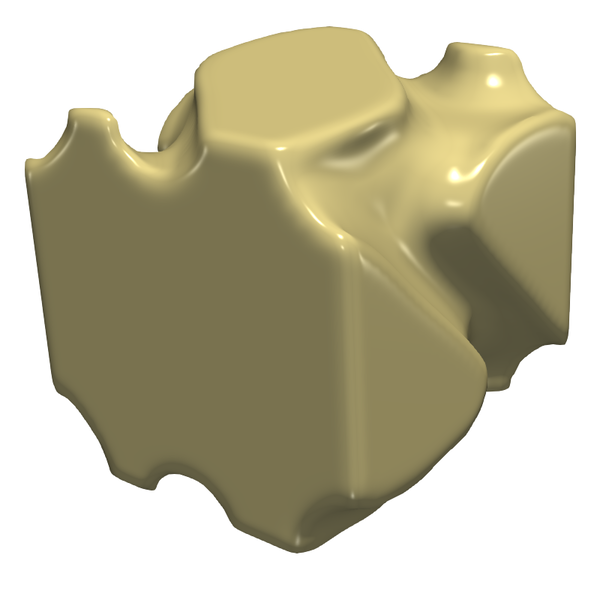}}  \end{minipage}
       \\
          \begin{minipage}{0.25\textwidth} {\includegraphics[width=0.9\textwidth]{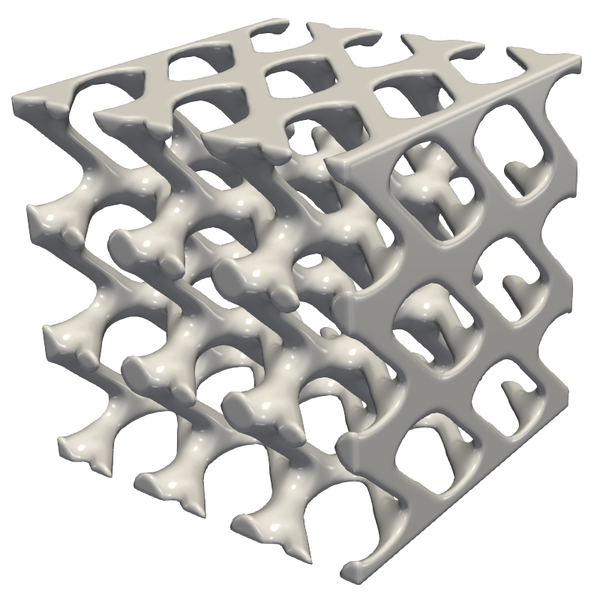}} \end{minipage}
       &  \begin{minipage}{0.25\textwidth} {\includegraphics[width=0.9\textwidth]{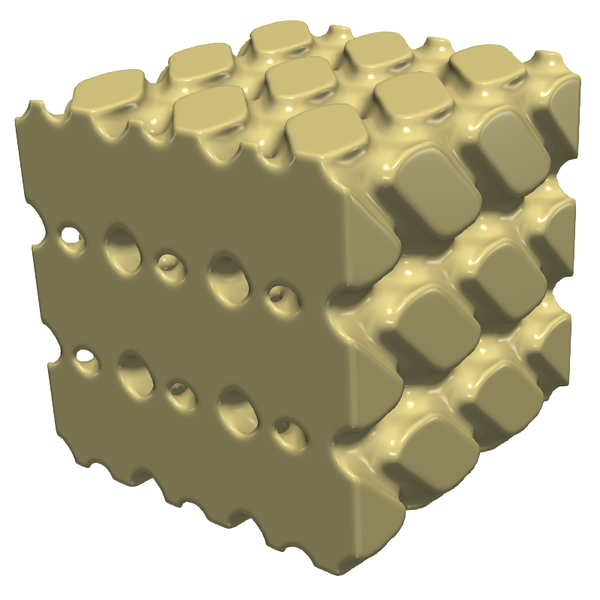}} \end{minipage}
       \\
  \end{tabular}
  
  \begin{tabular}{ c  c  c }
    \multicolumn{3}{c}{von Mises stresses bone} \\ 
      \begin{minipage}{0.25\textwidth} {\includegraphics[width=0.9\textwidth]{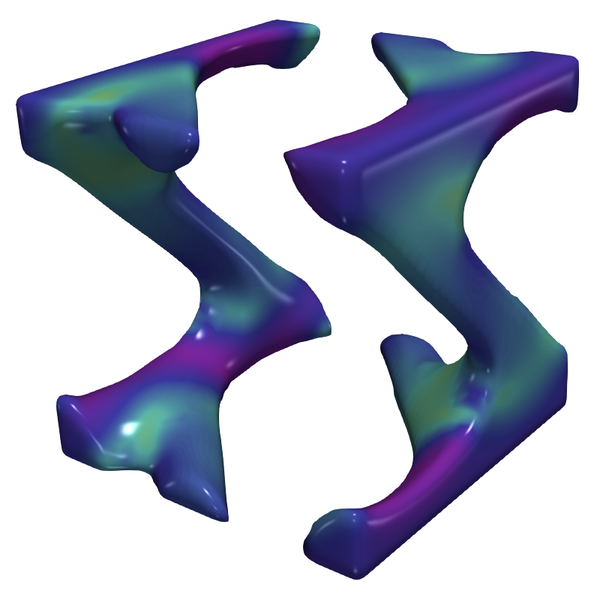}} \end{minipage}
    & \begin{minipage}{0.25\textwidth} {\includegraphics[width=0.9\textwidth]{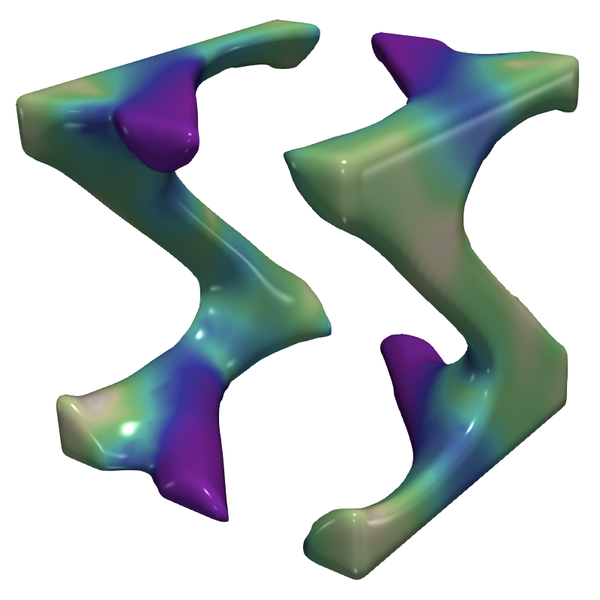}} \end{minipage}
    & \begin{minipage}{0.25\textwidth} {\includegraphics[width=0.9\textwidth]{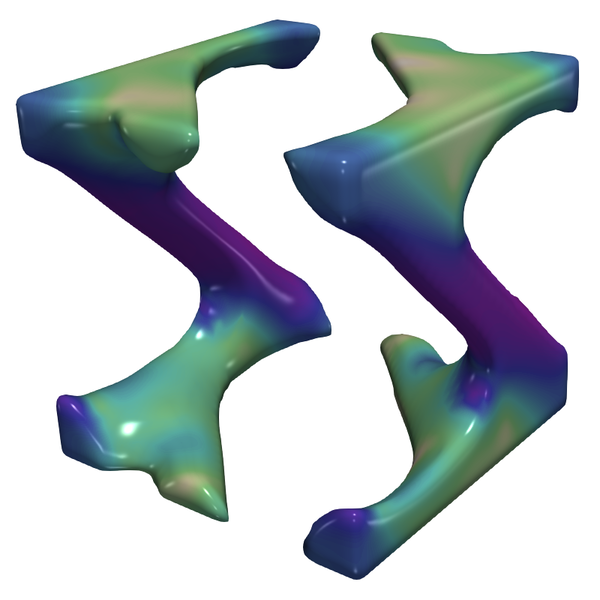}} \end{minipage}
    \\
    \multicolumn{3}{c}{von Mises stresses polymer} \\
      \begin{minipage}{0.25\textwidth} {\includegraphics[width=0.9\textwidth]{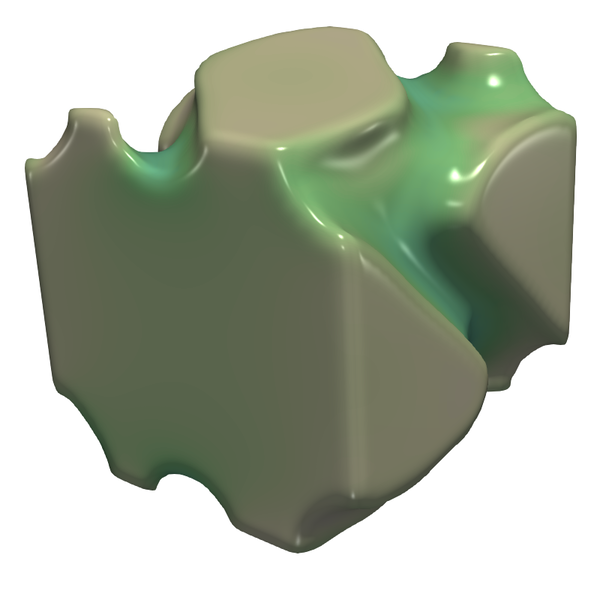}} \end{minipage}
    & \begin{minipage}{0.25\textwidth} {\includegraphics[width=0.9\textwidth]{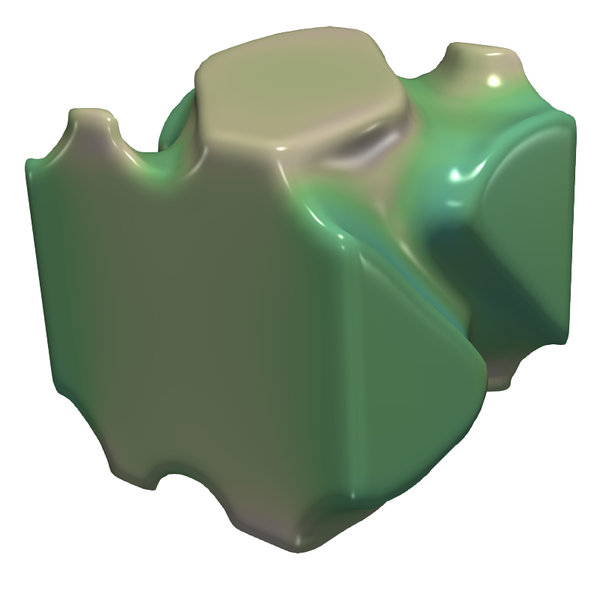}} \end{minipage}
    & \begin{minipage}{0.25\textwidth} {\includegraphics[width=0.9\textwidth]{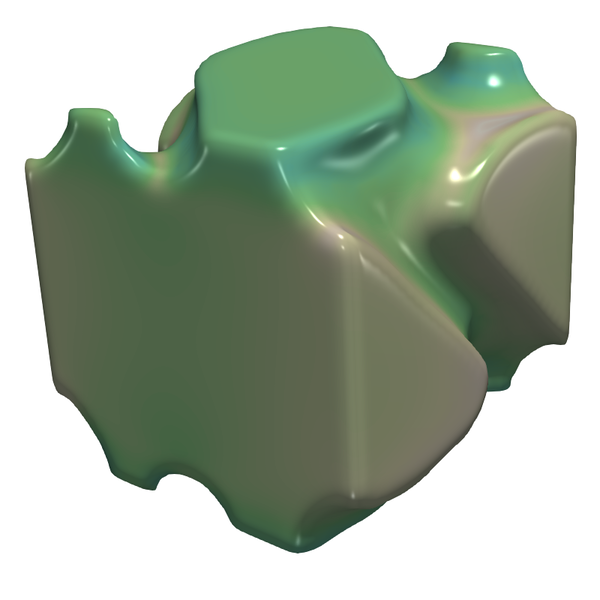}} \end{minipage}
    \\
  \end{tabular}
}
 \caption{\label{fig:bp} Optimal bone and polymer micro-structures with color coded von Mises stresses using a log scaled color value in HSV model.}
\end{figure}

\paragraph{Realistic material parameters for polymer and bone.}
Real bone is substantially stiffer than the bioresorbable polymer with a  $15$ times larger Youngs modulus and the Poisson ratio of  $\nu^B  = 0.1$ compared to the Poisson ratio  $\nu^P = 0.3$ for the polymer.
Figure \ref{fig:bp} shows the optimal bone and polymer subdomains together with a plot of the von Mises stresses on the boundary of the corresponding subdomains in the fundamental cell. 
Here, again the case of one compression load and two shear loads is taken into account.

\begin{figure}[h]
\begin{center}
\includegraphics[width=0.22\textwidth]{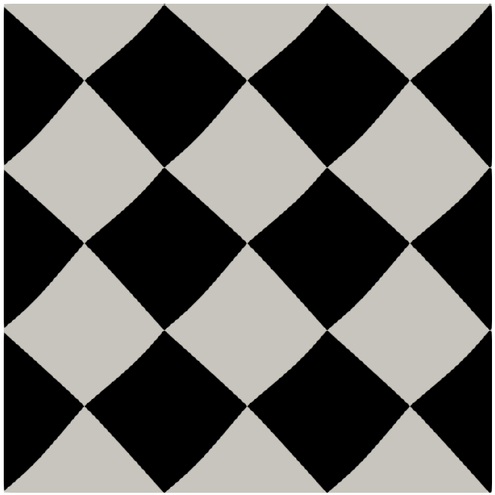}\hspace{1ex}
\includegraphics[width=0.22\textwidth]{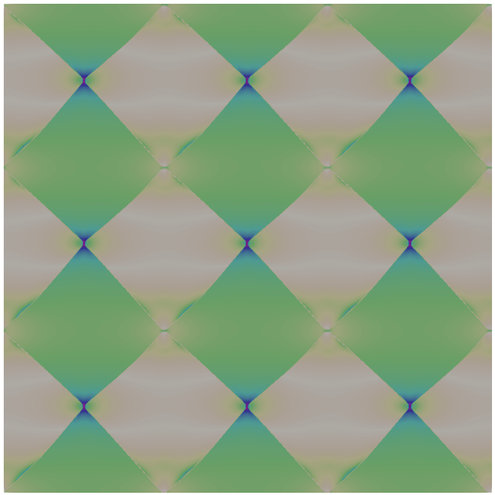}\hspace{1ex}
\includegraphics[width=0.22\textwidth]{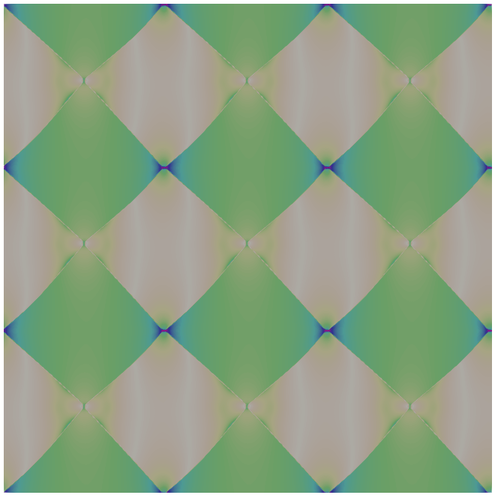}\hspace{1ex}
\end{center}
\caption{\label{fig:2D} An optimal domain decompositions in 2D in case of the hard--soft approximation with $\delta=10^{-4}$ is depicted  
for a load scenario with two different loads corresponding to two compression modes ($C^{1111}_*$ and $C^{2222}_*$).
A block of $3\,x\,3$ cells is plotted with the two subdomains in white and black together with a color plot of the von Mises stresses. 
}
\end{figure}
\begin{Remark}
 Let us briefly comment on the 2D case ($d=2$). Figure~\ref{fig:2D} shows the numerical result for a scenario with two uniaxial compression loads. 
 In the optimized shape configuration we obtain  diamond shaped regions of both subdomains which meet at the tips of the diamonds. 
In the context of our hard--soft approximation introduced in Section~\ref{sec:exist} this is a mechanically admissible configuration.
For a hard--void shape optimization model and two uniaxial compression loads in vertical and horizontal direction
no mechanically favourable splitting of the unit square $[0,1]^2$ into two subdomains is possibly.
 Indeed, a uniaxial load requires a truss with non vanishing interior connecting the components of the boundary opposite in the loading direction.
 A truss configuration simultaneously in horizontal and vertical direction for both subdomains is thus topologically impossible.

\end{Remark}

\section{Conclusions} \label{sec:conclusions}
We considered the problem of designing an optimal periodic microstructure for a domain splitting problem in shape optimization. The setting of this article is motivated by the biomechanical application of designing optimal scaffolds for bone regeneration, where both the scaffold as well as the regenerated bone (which can grow only in the space not occupied by the scaffold material) need to be stable individually. The numerical method presented here is able to treat a general cost functional whose input is given by the effective  moduli of an elastic material occupying either the optimized domain or its complement. 

The case of maximizing the compressive moduli in the three coordinate axes can be compared to the scalar problem of maximizing for example isotropic heat as well as electrical conduction in a two phase material where in each phase one parameter is large and the other is small \cite{Torquato:2002gh}.
Our simulations suggest that also in this 3d-elasticity setting a domain separation by a periodic minimal surface appears not to be optimal but with already comparable small cost, 
which was verified for the scalar case in \cite{Silvestre:2007hn}.

The optimization for other load cases, in particular a compression in one direction combined with a shear in the two directions orthogonal to the compressive load, 
yields optimal structures that are very clearly distinct from minimal surfaces.
Since this is the most physiologically relevant case, 
as the typical loading condition on major long bones is compression and torsion,
we note that it might be possible to further improve scaffold designs based on minimal surfaces (see Figure \ref{fig:PrintedScaffolds} and \cite{Kapfer:2011kz}), which are currently considered for medical practice.

A number of open issues remain. 
So far, we investigated solely spatially homogeneous microstructure and optimized the periodic scaffold on the fundamental cell. 
When considering realistic patient specific implant geometries and corresponding boundary conditions 
one has to exploit  the inherent multiscale nature of the problem and design an optimal microstructure with spatially varying microscopic shapes.
Furthermore, a major drawback of our current approach is that some significant physiological conditions are not considered in this work. 
A scaffold design as seen in Figure \ref{fig:bp} may be optimal from a purely mechanical point of view, but the very low porosity would seriously impede vascularization and therefore prevent the regeneration of bone matrix. 
Thus, in future study, quantities like the effective diffusivity in the pores should be considered (as it is done for example done in \cite{Adachi06}). 
The uncertainty in translating a given scaffold design into an additively manufactured scaffold can be observed in Figure \ref{fig:PrintedScaffolds}, where strands and layers stemming from the printing process are clearly visible.
A quantification of the uncertainty in the resulting effective elastic moduli due to both periodic as well as random fluctuations in the final product thus seems necessary as well.

\section*{Acknowledgements} We thank D.~Valainis (TU-Munich) for providing photographs of additively manufactured bone scaffold designs. P.~Dondl acknowledges support from the German Scholars Organization/Carl-Zeiss-Stiftung via the ``Wissenschaftler-R{\"u}ckkehrprogramm''.
M.~Rumpf and S.~Simon acknowledge support of the Hausdorff Center for Mathematics and the 
Collaborative Research Center 1060 funded by the German Research Foundation.

\bibliographystyle{alpha}
\bibliography{references,all,own}

\end{document}